\documentclass[12pt]{amsart}

\usepackage[margin=1.3in]{geometry}
\usepackage[utf8]{inputenc}
\usepackage[T1]{fontenc}
\usepackage{comment}
\usepackage{graphicx}
\usepackage{csquotes}
\usepackage{microtype}
\usepackage[skip=.5\baselineskip]{parskip}
\usepackage{anyfontsize}

\usepackage[dvipsnames]{xcolor}
\definecolor{highlight}{cmyk}{0.90,0.99,0,0}

\usepackage[textsize=footnotesize]{todonotes}
\setuptodonotes{color=blue!30}
\usepackage{amsmath,amssymb,amsthm,mathtools}
\usepackage{mathrsfs}
\usepackage{euscript}
\usepackage[cal=cm]{mathalfa}
\usepackage{xifthen}
\usepackage{thmtools}
\numberwithin{equation}{section}

\usepackage[lining]{libertine}
\usepackage{courier}
\usepackage[libertine,libaltvw,liby]{newtxmath}
\makeatletter
\renewcommand*\libertine@figurestyle{LF}
\makeatother

\usepackage{enumitem}
\setlist[enumerate,1]{label=(\roman*),itemsep=0.9ex}
\setlist[itemize]{itemsep=0.9ex}

\usepackage[colorlinks=true,allcolors=highlight,hyperfootnotes=false,linktocpage=true]{hyperref}
\usepackage{url}
\usepackage{cleveref}

\usepackage{tikz,tikz-cd}
\usepackage{pgfplots}
\pgfplotsset{compat=1.18}
\usetikzlibrary{shapes.misc}
\tikzset{
    main node/.style={circle,draw,inner sep=2pt,minimum size=0.8cm},
    crossed node/.style={cross out,draw,minimum size=6pt,inner sep=0pt},
    caption node/.style={font=\large},
    blue label/.style={blue,font=\footnotesize}}

\hypersetup{bookmarksdepth = 3}

\makeatletter
\def\@tocline#1#2#3#4#5#6#7{\relax
  \ifnum #1>\c@tocdepth
  \else
    \par \addpenalty\@secpenalty\addvspace{#2}%
    \begingroup \hyphenpenalty\@M
    \@ifempty{#4}{%
      \@tempdima\csname r@tocindent\number#1\endcsname\relax
    }{%
      \@tempdima#4\relax
    }%
    \parindent\z@ \leftskip#3\relax \advance\leftskip\@tempdima\relax
    \rightskip\@pnumwidth plus4em \parfillskip-\@pnumwidth
    #5\leavevmode\hskip-\@tempdima
    \ifcase #1
    \or\or \hskip 1em \or \hskip 2em \else \hskip 3em \fi
    #6\nobreak\relax
    \dotfill\hbox to\@pnumwidth{\@tocpagenum{#7}}\par
    \nobreak
    \endgroup
  \fi}
\makeatother

\DeclareMathOperator{\Tot}{Tot}

\newcommand{\nRootSpec}{cyclotomic specialisation}
\newcommand{\eee}{\epsilon}
\newcommand{\ee}{\epsilon}
\newcommand{\ddd}{\delta}
\newcommand{\dd}{\delta}

\newcommand{\mr}{\mathrm}

\newcommand{\re}{\mathrm{e}}

\newcommand{\modulifont}[1]{#1}
\newcommand{\Mbar}{\overline{\modulifont{M}}}
\newcommand{\vir}{\mathrm{virt}}

\newcommand{\Gm}[1][]{%
  \ifthenelse{\isempty{#1}}
  {\mathbb{C}^\times}
  {(\mathbb{C}^\times)^{#1}}
}

\makeatletter
\def\makeCal#1{\expandafter\newcommand\csname c#1\endcsname{\mathcal{#1}}}
\def\makeBB#1{\expandafter\newcommand\csname b#1\endcsname{\mathbb{#1}}}
\def\makeFrak#1{\expandafter\newcommand\csname f#1\endcsname{\mathfrak{#1}}}
\def\makeScr#1{\expandafter\newcommand\csname s#1\endcsname{\mathscr{#1}}}
\def\makeSF#1{\expandafter\newcommand\csname sf#1\endcsname{\mathsf{#1}}}

\count@=0
\loop
  \advance\count@ 1
  \edef\y{\@Alph\count@}
  \expandafter\makeCal\y
  \expandafter\makeBB\y
  \expandafter\makeFrak\y
  \expandafter\makeScr\y
  \expandafter\makeSF\y
\ifnum\count@<26
\repeat
\makeatother

\theoremstyle{plain}
\newtheorem{thm}{Theorem}[section]
\newtheorem{prop}[thm]{Proposition}
\newtheorem{lem}[thm]{Lemma}

\theoremstyle{definition}

\newtheorem{rmk}[thm]{Remark}
\newtheorem{example}[thm]{Example}

\theoremstyle{plain}
\newtheorem{innercustomthm}{Theorem}
\newenvironment{customthm}[1]
{\renewcommand\theinnercustomthm{#1}\innercustomthm}
{\endinnercustomthm}

\newtheorem{innercustomcor}{Corollary}
\newenvironment{customcor}[1]
{\renewcommand\theinnercustomcor{#1}\innercustomcor}
{\endinnercustomcor}

\newtheorem{innercustomconj}{Conjecture}
\newenvironment{customconj}[1]
{\renewcommand\theinnercustomconj{#1}\innercustomconj}
{\endinnercustomconj}

\AddToHook{cmd/appendix/before}{\crefalias{section}{appendix}}

\crefname{equation}{Equation}{Equations}
\crefname{eqnarray}{Equation}{Equations}

\crefname{thm}{Theorem}{Theorems}
\crefname{prop}{Proposition}{Propositions}
\crefname{lem}{Lemma}{Lemmas}
\crefname{cor}{Corollary}{Corollaries}
\crefname{conj}{Conjecture}{Conjectures}
\crefname{defn}{Definition}{Definitions}
\crefname{rmk}{Remark}{Remarks}
\crefname{example}{Example}{Examples}
\crefname{defprop}{Definition/Proposition}{Definitions/Propositions}

\crefname{innercustomthm}{Theorem}{Theorems}
\crefname{innercustomcor}{Corollary}{Corollaries}
\crefname{innercustomconj}{Conjecture}{Conjectures}

\crefname{section}{Section}{Sections}
\crefname{appendix}{Appendix}{Appendices}
\crefname{figure}{Figure}{Figures}
\crefname{table}{Table}{Tables}

\usepackage[
  backend=biber,
  style=alphabetic,
  natbib=true,
  url=false,
  doi=false,
  eprint=true,
  isbn=false,
  maxcitenames=5,
  maxbibnames=4,
  maxalphanames=4,
  sortcites=true
]{biblatex}

\AtEveryBibitem{\clearlist{language}}
\renewbibmacro{in:}{%
  \ifboolexpr{
    test {\ifentrytype{article}}
    or
    test {\ifentrytype{inproceedings}}
  }
  {}
  {\printtext{\bibstring{in}\intitlepunct}}%
}

\DeclareFieldFormat[article,inbook,incollection,inproceedings,patent,thesis,unpublished,techreport,misc,book]{title}{\mkbibquote{#1}}

\title[Gromov--Witten theory of the local model of the standard flop]{Gromov--Witten theory of the local model\\ of the standard flop}

\author[A. D\"{a}nzer]{Alex D\"{a}nzer}
\address{Department Mathematik, ETH Z\"urich, CH-8006 Z\"urich, Switzerland}
\email[AD]{adaenzer@student.ethz.ch}

\author[Y. Schuler]{Yannik Schuler}
\email[YS]{yannik.schuler@math.ethz.ch}

\begin{document}

\begin{abstract}
    We conjecture that the generating series of equivariant Gromov--Witten invariants of the local Calabi--Yau $(2n-1)$-fold $\operatorname{Tot} \mathcal{O}_{\mathbb{P}^{n-1}} (-1)^{\oplus n}$ at fixed genus can be presented as a linear combination of polylogarithms, generalising a well-known property of the resolved conifold. We prove this conjecture for two specialisations of the equivariant parameters. In the first, the parameters are specialised to roots of unity and the proof proceeds via Givental--Teleman reconstruction. In the second, we obtain a closed-form evaluation of the generating series from which the conjecture follows. For $n=3$, the resulting closed formula establishes a special case of a Gopakumar--Vafa-type conjecture for Calabi--Yau fivefolds.
\end{abstract}

\maketitle
\tableofcontents

\section*{Introduction}
This note concerns the Gromov--Witten theory of the local Calabi--Yau $(2n-1)$-fold
\begin{equation*}
    Z_n = \Tot \mathcal{O}_{\mathbb{P}^{n-1}}(-1)^{\oplus n},\qquad n>1.
\end{equation*}
This variety appears in the literature as the local model for the standard flop of Bondal and Orlov \cite{BO95:SemiorthDec} and is a higher-dimensional generalisation of the resolved conifold. We equip $Z_n$ with the action of the torus $T = \Gm[n] \times \Gm[n]$. We let the first factor $\Gm[n]$ of the torus act on the homogeneous coordinates of the base $\mathbb{P}^{n-1}$ with weights $\eee_0,\ldots,\eee_{n-1}$ and we consider the natural lift of this action to $Z_n$. The second factor of the torus acts fibre-wise with weights $\ddd_0,\ldots,\ddd_{n-1}$ on the respective line bundle factors.

We consider the all-degree fixed-genus $g$ generating series of $T$-equivariant Gromov--Witten invariants of $Z_n$:
\begin{equation}
    \label{eq: F g intro}
    F^{Z_n, T}_g(Q) \coloneqq \sum_{d = 1}^{\infty} Q^d  \int_{[\Mbar_g(Z_n,d)]^{\vir}_T}1.
\end{equation}
A priori, this is just a formal power series in $Q$ whose coefficients are polynomials in $\eee_0,\ldots, \eee_{n-1}$ and $\ddd_{0},\ldots,\ddd_{n-1}$, but we conjecture the series can be presented as a linear combination of polylogarithms
\begin{equation*}
    \mr{Li}_{-k}(Q) = Q + 2^k Q^2 + 3^k Q^3 + \ldots .
\end{equation*}

\begin{customconj}{A}
    \label{conj: F GW via Li}
    For all $g,m\geq 0$ with $2g-2+m>0$ there exits some $N_{g,n}>0$ such that the series $$(Q \partial_Q)^m F^{Z_n, T}_g$$ lies in the linear $\mathbb{Q}[\eee_0,\ldots, \eee_{n-1},\ddd_{0},\ldots,\ddd_{n-1}]$-span of the polylogarithms
    \begin{equation*}
        \mr{Li}_{0}\big((-1)^nQ\big), \quad \mr{Li}_{-1}\big((-1)^nQ\big), \quad \ldots, \quad \mr{Li}_{-N_{g,n}-m}\big((-1)^nQ\big).
    \end{equation*}
    Furthermore, the series is symmetric under permutation of $\eee_0,\ldots, \eee_{n-1}$ and $\ddd_{0},\ldots,\ddd_{n-1}$ and, after analytic continuation, it has parity $(-1)^m$ under the simultaneous substitution
    \begin{equation}
        \label{eq: Q to Q sym intro}
        Q \longmapsto Q^{-1}, \qquad (\eee_0,\ldots, \eee_{n-1}) \longleftrightarrow (\ddd_{0},\ldots,\ddd_{n-1}).
    \end{equation}
\end{customconj}


We prove our conjecture in the following two cases.

\begin{customthm}{B}
    \label{thm: F GW via Li}
    \Cref{conj: F GW via Li} holds under the following weight specialisations:
    \begin{enumerate}[leftmargin=2em]
        \item\label{item: thm B i} in the {\nRootSpec} $\eee_i = \zeta^{i} \ee$ and $\ddd_i = \zeta^{i} \dd$ with $\zeta = \mr{e}^{2\pi \mr{i}/n}$, $i\in\{0,\ldots, n-1\}$, the conjecture holds with the bound
        \begin{equation*}
            N_{g,n} = \left\lfloor \frac{(n-1)(4g-4)}{n}\right\rfloor -1.
        \end{equation*}
        
        \item\label{item: thm B ii} in the {anti-diagonal specialisation} $\eee_i = -\ddd_i$, $i\in\{0,\ldots, n-1\}$ it holds with
        \begin{equation*}
            N_{g,n} = 2g-3.
        \end{equation*}
    \end{enumerate}
\end{customthm}

\subsection{The \nRootSpec}
We prove \Cref{conj: F GW via Li} in the \nRootSpec\ \ref{item: thm B i} via Givental--Teleman reconstruction \cite{giventalconjecture,Telemanproof}. The proof makes up the majority of this paper and is spread over Sections \ref{sec: setup}--\ref{sec: GivTel}. We start with a careful analysis of the small quantum cohomology in the \nRootSpec. We determine most structure constants of the quantum product and constrain the $Q$-dependence of an asymptotic expansion of the small $J$-function via its Picard--Fuchs equation (\Cref{sec: QuCoh}). This is the place where the \nRootSpec\ is crucially used in our proof as it greatly simplifies the Picard--Fuchs operator. This step is therefore the main obstruction towards generalising our proof of \Cref{thm: F GW via Li}\ref{item: thm B i} to arbitrary torus actions, cf.~\Cref{rmk: obstruction generalisation}. The same difficulty of going beyond \nRootSpec s is also present in the work of Lho--Pandharipande and Lho \cite{LhoLocalP1,CrepantresolutionC3/Z3,LP18:KP2,Lho21:KPn} and Genlik--Tseng \cite{GT24:Cn/ZnI,GT26:Cn/ZnII} concerning the Gromov--Witten theory of $K_{\bP^{n-1}}$ and $[\bC^{n}/\mathbb{Z}_{n}]$. It should be remarked that the approach we pursue in this paper to prove \Cref{thm: F GW via Li}\ref{item: thm B i} very much follows the ideas in loc.~cit.

Our analysis of the small $J$-function is then used to constrain the $Q$-dependence of the $R$-matrix of $Z_n$ akin to the statement of \Cref{conj: F GW via Li}. This result is then transported to the generating series \smash{$F_{\mathrlap{g}}{\mathstrut}^{Z_n, T}$} via Givental--Teleman reconstruction, concluding our proof of \Cref{thm: F GW via Li}\ref{item: thm B i}. Details are found in \Cref{sec: GivTel}.

We remark that the symmetry $Q\mapsto Q^{-1}$, together with the simultaneous exchange of the equivariant parameters $(\eee_i)\leftrightarrow(\ddd_i)$, is a higher-genus manifestation of the crepant transformation property of Gromov--Witten series \cite{Rua02:CRC,BG09:CRC,CR13:CRC}. It generalises the well-known flop invariance of the resolved conifold, $n=2$, to higher $n$. To recall the construction of the flop, note that a point of $Z_n$ is a line $\ell\subset\bC^n$ together with vectors $v_1,\ldots,v_n\in\ell$. Forgetting the line and arranging these vectors as the columns of an $n\times n$ matrix defines a morphism $Z_n\to\bC^{n\times n}$ whose image is the subvariety $Y_n$ of matrices of rank at most one. The variety $Y_n$ is singular at the origin and $Z_n\to Y_n$ is a crepant resolution. Exchanging the roles of columns and rows gives a second resolution $Z_n'\to Y_n$ with $Z_n\cong Z_n'$ and the induced birational map is the so-called standard \mbox{$\bP^{n-1}$-flop}
\begin{equation*}
    \begin{tikzcd}
        Z_n \ar[rr,dashed] \ar[dr] & & Z_n' \ar[dl]\\
        & Y_n &
    \end{tikzcd}
\end{equation*}
Since the flop exchanges the roles of the equivariant parameters $\eee_i$ and $\ddd_i$ of the torus action, this provides a conceptual explanation for the self-symmetry \eqref{eq: Q to Q sym intro} of the generating series of Gromov--Witten invariants. The invariance of quantum cohomology under standard $\bP^{n-1}$-flops, and generalisations thereof, was previously investigated by Lee--Lin--Wang \cite{LLW10:FlopsQuRings,LLW16:FlopInvarI,LLW16:FlopInvarII} and Lee--Lin--Qu--Wang \cite{LLQW16:FlopInvarIII}.

\subsection{The anti-diagonal specialisation}
In the anti-diagonal specialisation \ref{item: thm B ii}, \Cref{thm: F GW via Li} is proven via a direct evaluation of the series.

\begin{customthm}{C}
    \label{thm: closed form}
    In the anti-diagonal specialisation, for fixed degree $d>0$ we have
    \begin{equation*}
        \sum_{g \geq 0} u^{2g-2} \int_{[\Mbar_g(Z_n,d)]^{\vir}_T} 1 = \frac{(-1)^{n(d+1)}}{d} \sum_{0\leq i<j <n} ~ \frac{1}{\prod_{k\in\{i,j\}} 2 \sin \big(\frac{du}{2} \prod_{l\neq i,j} (\eee_l - \eee_k)\big)}
    \end{equation*}
    when $n$ is even. A formula valid for odd $n$ is obtained by substituting the above sine function by a hyperbolic sine, cf.~equation \eqref{eq: closed form n odd}.
\end{customthm}

We prove this formula in \Cref{sec: closed form}. \Cref{thm: F GW via Li}\ref{item: thm B ii} then follows from the observation that the $u^{2g-2}$-coefficient of the right-hand side of the formula is a $d$-independent polynomial in $\eee_0,\ldots,\eee_{n-1}$ times a factor $(-1)^{nd}d^{2g-3}$. This means the fixed-genus generating series must be proportional to $\mr{Li}_{3-2g}((-1)^n Q)$. In particular, \Cref{conj: F GW via Li} holds in the anti-diagonal specialisation. For $n=2$ this reproduces the well-known formula for the resolved conifold
\begin{equation*}
    F^{Z_2, T}_g(Q) = \frac{|B_{2g}|}{2g \, (2g-2)!}\, \mr{Li}_{3-2g}(Q),\qquad g \geq 2
\end{equation*}
proven by Faber and Pandharipande \cite[Thm.~3]{FP00:HodgeIntGW}. Indeed, our proof of \Cref{thm: closed form} is a direct generalisation of Faber and Pandharipande's proof to $n>2$.

Let us now discuss the special case $n=3$. Taking the limit $\eee_i\rightarrow 0$, \Cref{thm: closed form} reproduces a formula proven by Pandharipande--Zinger for the genus-one Gromov--Witten invariants of $Z_3$ \cite[Lem.~3.1]{PandharipandeNumerative5folds}:
\begin{equation}
    \label{eq: PZ formula}
    \int_{[\Mbar_1(Z_3,d)]^{\vir}_T} 1 = \frac{(-1)^d}{8d}.
\end{equation}
A similar formula for arbitrary $n$ was proven by Hu in \cite[Sec.~3.0.1]{Hu17:g1GWlocPn}. Beyond genus one, the closed formula for $Z_3$ in \Cref{thm: closed form} establishes an instance of a conjecture of the second author \cite[Conj.~A]{Sch26:TVCY5} concerning the equivariant Gromov--Witten theory of Calabi--Yau fivefolds. Specialised to $Z_3$, the conjecture asserts the existence of rational functions
\begin{equation*}
    \Omega_d(q_0,\ldots,q_{5}) \in \mathbb{Z}\big[\tfrac{1}{2}\big]\!\left[q_0^{\pm1/2},\ldots, q_{5}^{\pm1/2}, \left\{ \frac{1}{1-\prod_i q_i^{a_i/2}}\right\}_{\boldsymbol{a}\in \mathbb{Z}^{6} \setminus\{0\}}\right]
\end{equation*}
such that under the change of variables $q_i = \mr{e}^{u \eee_i}$ and $q_{i+3} = \mr{e}^{u \ddd_i}$, $i\in\{0,1, 2\}$:
\begin{equation}
    \label{eq: num GW M2}
    \sum_{g \geq 0} u^{2g-2} \int_{[\Mbar_g(Z_3,d)]^{\vir}_T} 1 = \sum_{k \mid d} \frac{1}{k} \, \Omega_{d/k}(q_i^k).
\end{equation}
\Cref{thm: closed form} establishes this conjecture with an explicit expression for $\Omega_d$.

\begin{customcor}{D}
    Conjecture A of \cite{Sch26:TVCY5}, i.e. equation \eqref{eq: num GW M2}, holds for $Z_3$ in the anti-diagonal specialisation \ref{item: thm B ii} with
    \begin{equation*}
        \Omega_1 = -\frac{\big(q_0 q_1 q_2\big)^{1/2}}{\big(q_0^{1/2}+q_1^{1/2}\big)\big(q_1^{1/2}+q_2^{1/2}\big)\big(q_2^{1/2}+q_0^{1/2}\big)} ,\qquad \Omega_2 = \frac{q_0 q_1 q_2}{\big(q_0+q_1\big)\big(q_1+q_2\big)\big(q_2+q_0\big)}
    \end{equation*}
    and $\Omega_{d} = 0$ for $d\geq 3$.
\end{customcor}

\subsection{Beyond the two specialisations}
Holmes and the second author investigated the Gromov--Witten invariants of $Z_3$ in the Calabi--Yau specialisation of the equivariant parameters $\sum_i \ee_i + \sum_i \dd_i = 0$. The following formula for the genus-$2$ series was predicted based on numerical calculations \cite[Sec.~4.1.3]{HS26:ExpMMS}:
\begin{equation*}
    \begin{split}
        & \qquad F^{Z_3, T}_2(Q) = \\
        & \frac{1}{64} \left(3\!\!\!\sum_{0\leq i<j\leq 2}\!\!\! (\eee_i \eee_j-\ddd_i\ddd_j) + 2 (\eee_0+\eee_1+\eee_2)^2\right)\,\mr{Li}_{-1}(-Q) + \frac{3}{64} \left(\sum_{0\leq i<j\leq 2} \!\!\!\eee_i \eee_j-\ddd_i\ddd_j\right)\,\mr{Li}_{-2}(-Q)\,. 
    \end{split}
\end{equation*}
There are two lessons to learn from this example: First, one should generally expect a linear combination of several polylogarithms and not just a single one as in the anti-diagonal specialisation \ref{item: thm B ii}. Second, the bound $N_{g,n}$ in \Cref{thm: F GW via Li}\ref{item: thm B ii} does not hold in general as it is off by one in the above example. We leave it as an open problem for future research to determine a sharp bound $N_{g,n}$ valid for generic equivariant parameters.

\subsection*{AI disclosure}
The mathematical content of the paper is due to the authors. Generative AI was used to test the $L$-degree bound and streamline the proof of \Cref{lemmwwithlogterm} and to debug code used for the computations in \Cref{example:n3}. Apart from this, generative AI was solely used to improve the language and grammar of the paper.

\subsection*{Acknowledgements}
We thank Andrea Brini, Alessandro Giacchetto and Rahul Pandharipande for fruitful discussions and helpful comments. 
The second author was supported by SNF-200020-219369 and a Walter Benjamin Fellowship of the Deutsche For\-schungs\-ge\-mein\-schaft (DFG) --- Projektnummer 576663726.

\section{Setup and notation}
\label{sec: setup}

\subsection{Geometric preliminaries}
The following exposition closely follows Section~5 of \cite{BS24:refGW}. Recall from the introduction that for $n > 1$ we denote by $Z_n$ the total space $\operatorname{Tot} (\mathcal{O}_{\mathbb{P}^{n-1}}(-1)^{\oplus n})$. We can present $Z_n$ as a GIT quotient of $\mathbb{C}^{2n}$ with coordinates $(x_0,\ldots,x_{n-1};y_0,\ldots,y_{n-1})$:
\begin{equation}
    Z_n = (\mathbb{C}^{2n} \backslash \Delta)\big/\mathbb{C}^\times.
    \nonumber
\end{equation}
Here, $\Delta \coloneqq \{(x,y) \in \mathbb{C}^{2n} \mid x_0=\ldots=x_{n-1}=0\}$ is the unstable locus and the torus $\mathbb{C}^{\times}$ acts on $\mathbb{C}^{2n}$ via
\begin{align}
    \mathbb{C}^{\times}\times \mathbb{C}^{2n} &\longrightarrow \mathbb{C}^{2n}\,, \notag\\
    (\lambda, x_0,\ldots, x_{n-1}, y_0,\ldots, y_{n-1}) &\longmapsto 
    (\lambda x_0, \ldots, \lambda x_{n-1}, \lambda^{-1} y_0,\ldots,\lambda^{-1} y_{n-1})\,.
    \nonumber
\end{align}
The first $n$ coordinates $[x_0:\ldots:x_{n-1}]$ correspond to the homogeneous coordinates of the base $\mathbb{P}^{n-1}$, while $(y_0,\ldots, y_{n-1})$ are the coordinates of the fibre $Z_n \rightarrow \mathbb{P}^{n-1}$. The natural action of the torus $T \simeq (\mathbb{C}^{\times})^{2n}$ on $\mathbb{C}^{2n}$, whose tangent weights we denote by
\begin{equation}
    (\eee_0,\ldots,\eee_{n-1},\ddd_0,\ldots,\ddd_{n-1}),
    \nonumber
\end{equation}
descends to an action on the quotients $Z_n$ and $\mathbb{P}^{n-1}$. The action has $n$ isolated fixed points, namely the coordinate points
\[ P_i = [0:\ldots:0:1:0:\ldots:0]
\]
of the projective space with $1$ at the $i$th position, $i\in \{0,\ldots,n-1\}$.

The torus action has a natural lift to the line bundle $\mathcal{O}_{\mathbb{P}^{n-1}}(1)$, where the torus acts on the fibre of each fixed point $P_i$ with weight $-\eee_i$, $i \in \{0,\ldots,n-1\}$. 
We denote by $H \coloneqq \pi^{*} c^{T}_{1}\!\big(\mathcal{O}_{\mathbb{P}^{n-1}}(1)\big)$ the associated equivariant hyperplane class. Then the $T$-equivariant Chow cohomology of $Z_n$, which we take with complex coefficients, is given by
\begin{equation}
    \mathrm{CH}_{T}^{\bullet}(Z_n)\, \cong \, 
    \mathbb{C}[H, \eee_0,\ldots,\eee_{n-1},\ddd_{0},\ldots,\ddd_{n-1}] \, \Big/ \,
    \Big( {\textstyle\prod_{i=0}^{n-1}} (H + \eee_i)\Big).
    \nonumber
\end{equation}
One computes the $T$-characters on the tangent spaces $T_{P_i}Z_n$ to be
\begin{equation}
    w_{i}^j =
    \begin{cases}
    \eee_{(1+i+j)\bmod n} - \eee_i\,, & j \in \{0,1,\ldots,n-2\}, \\[6pt]
    \eee_i + \ddd_{j-n+1}, & j \in \{n-1,\ldots,2(n-1)\},
    \end{cases}
    \label{TMcharacters}
\end{equation}
where $i\in \{0,\ldots,n-1\}$ and $j\in \{0,\ldots,2(n-1)\}$. Writing $\mathrm{w}_i \coloneqq e_T(T_{P_i}Z_n) = \prod_{j=0}^{2(n-1)}w_i^j$ for the Euler class of the tangent spaces, the equivariant intersection pairing on $Z_n$ is determined by the localisation formula:
\begin{equation}
    \eta(\phi,\psi) \coloneqq  \int_{[Z_n]_{T}}\phi  \cup \psi =\sum_{k=0}^{n-1}\frac{\phi|_{P_k}\cdot\psi|_{P_k}}{\mathrm{w}_k}.
    \label{eta def}
\end{equation}

\subsection{The \nRootSpec}
The subsequent sections \ref{sec: QuCoh} and \ref{sec: GivTel} will be concerned with the \emph{\nRootSpec} of the equivariant parameters:
\begin{equation}
\eee_i = \zeta^i \ee ,\qquad \ddd_i = \zeta^i \dd, \qquad \zeta \coloneqq \mr{e}^{\frac{2 \pi \mr{i}}{n}}.
\label{specialization}
\end{equation}
In this specialisation, the tangent Euler classes satisfy the relation $\mathrm{w}_i =\zeta^{(n-1)i}\mathrm{w}_0,$ greatly simplifying the sum in \cref{eta def}.
With respect to the (flat) basis \(\{H^a\}_{a\in\{0,\ldots,n-1\}}\) of $\mathrm{CH}_{T}^\bullet(Z_n)$, the entries of the Gram matrix of the intersection pairing are given by:
\begin{equation}
    \eta_{ab} = \frac{1}{\ee^n - (-1)^{n}\dd^n} \, \delta_{n-1,a+b}.
    \label{int_prod}
\end{equation}
We denote the components of the inverse matrix by $\eta^{ab} \coloneqq (\eta^{-1})_{ab}$.

\subsection{Gromov--Witten invariants} 
We work in the flat basis $\{H^a\}_{a\in \{0,\dots,n-1\}}$ generated by the equivariant hyperplane class $H$, and pick an element $\tau = \sum_{a=0}^{n-1}\tau_a H^a$ in the localised equivariant Chow ring of $Z_n$. For non-negative integers $a_1,\dots,a_m,r_1,\dots,r_m$ we introduce the correlator notation
\begin{align}
    \left\langle 
    H^{a_1}\psi_1^{r_1},\ldots, H^{a_m}\psi_m^{r_m} 
    \right\rangle^{Z_n}_{g,m,d}
    & \coloneqq
    \int_{[\Mbar_{g,m}(Z_n,d) ]^{\mathrm{vir}}_{T}}
    \prod_{i=1}^m \mathrm{ev}_i^*(H^{a_i})\,\psi_i^{r_i},
    \label{eq: correlator} \\[1em]
    \big\langle \!\!\big\langle
    H^{a_1}\psi_1^{r_1},\dots, H^{a_m}\psi_m^{r_m} 
    \big\rangle \!\! \big\rangle^{Z_n}_{g,m}
    & \coloneqq \!\!
    \sum_{\substack{d\geq 0 \\ k\geq 0}} \frac{1}{k!}
    \langle 
    H^{a_1}\psi_1^{r_1},\ldots, H^{a_m}\psi_m^{r_m},
    \underbrace{\tau,\ldots,\tau}_{k\ \text{times}}
    \rangle^{Z_n}_{g,m+k,d}.
    \nonumber
\end{align}
Later the following feature of these correlators will become crucial.

\begin{lem}
    \label{lem: s t dependence}
    In the \nRootSpec\ \eqref{specialization}, all correlators \eqref{eq: correlator}
    are rational functions in $\ee^n$ and $\dd^n$ and are polynomial in these variables if $d>0$.
\end{lem}

\begin{proof}
    Polynomiality in $\ee$ and $\dd$ is due to the properness of the moduli space of non-constant stable maps to $Z_n$. Moreover, all equivariant Gromov--Witten invariants must be invariant under the discrete symmetry group $\mathfrak{S}_{n}\times \mathfrak{S}_{n}$ acting on $Z_n$ by permuting the $n$ fixed points of the zero section $\mathbb{P}^{n-1}$ and permuting the $\mathcal{O}_{\mathbb{P}^{n-1}}(-1)$-factors of the vector bundle $Z_n$. As a consequence, all Gromov--Witten invariants can be expressed in elementary symmetric polynomials in $\eee_0,\ldots,\eee_{n-1}$ and $\ddd_0,\ldots, \ddd_{n-1}$. In the \nRootSpec\ all of them vanish apart from the top ones: $\ee^n$ and $\dd^n$.
\end{proof}

\subsection{Quantum product}
In the flat basis, the big quantum product is defined as
\begin{equation}
    \eta(H^{a} \circ_{\tau} H^{b} ,\, H^{c})
    \coloneqq \bigl\langle \!\!\bigl\langle H^{a},\, H^{b},\, H^{c} \bigr\rangle \!\! \bigr\rangle^{Z_n}_{0,3}
    \label{big quantum product}
\end{equation}
We denote its structure constants by $C^c_{ab}$:
\begin{equation}
     H^a \circ_{\tau} H^b = \sum_{c=0}^{n-1} C^c_{ab}(\tau) \, H^c.
     \nonumber
\end{equation}
Since the torus $T$ acts on $Z_n$ with finitely many fixed points, the quantum product is semi-simple. Let us, hence, write $\{\mathbf{e}_i\}_{i \in \{0, \dots, n-1\}}$ for the idempotents of the algebra, and $\mathbf{f}_i$ for 
their normalised counterparts:
\begin{equation}
    \mathbf{e}_i \circ_{\tau} \mathbf{e}_j = \delta_{i,j} \mathbf{e}_i, \quad \mathbf{f}_i = \Delta_i^{1/2}\mathbf{e}_i, \quad \text{where} \quad \Delta_i^{-1} = \eta(\mathbf{e}_i,\mathbf{e}_i).
    \nonumber
\end{equation}

\subsection{\texorpdfstring{$J$}{J}-function}
\label{sec: J fct}
We write
\begin{equation}
    J(\tau,z) \coloneqq (z + \tau_0) \mathbf{1} +\sum_{a=1}^{n-1}\tau_aH^a 
    + \sum_{a,b=0}^{n-1} H^a \, \eta^{ab} \, \Big\langle \! \!\Big\langle \frac{H^b}{z - \psi} \Big\rangle \!\! \Big\rangle_{0,1}^{Z_n}.
    \label{J-function}
\end{equation} 
for the big $J$-function of $Z_n$. Recall that
 \begin{equation}
     S(\tau,z) \coloneqq dJ(\tau,z)
     \nonumber
 \end{equation}
is then the fundamental solution to the Dubrovin connection. Its action on the flat basis is given by
\begin{equation}
    S(\tau,z)(H^a) = \partial_{\tau_a} J(\tau,z) = H^a + \sum_{b,c=0}^{n-1} H^c \, \eta^{cb} \, \Big\langle\!\!\Big\langle H^a,\frac{H^b}{z-\psi} \Big\rangle\!\!\Big\rangle_{0,2}^{Z_n}.
    \label{big S-funct}
\end{equation}

In the following, we will regularly restrict some of the above defined quantities to the small quantum locus, that is, the subspace where all deformation parameters $\tau_a$ vanish except for the one corresponding to the divisor class. In other words, one restricts to $\tau_a=0$ for all $a\neq1$ and so $\tau = \tau_1 H$. We adopt the convention that we use non-italic letters whenever we restrict to the small quantum cohomology locus, while we reserve italic letters for quantities defined on the big quantum cohomology locus. For instance, the restriction of the $J$-function \cref{J-function} to the small cohomology locus yields
\begin{equation}
     \mathrm{J}(Q,z) = zQ^{H/z}\left(\mathbf{1} + \sum_{d>0}\sum_{a,b} Q^dH^a\eta^{ab}\Big\langle \frac{H^b}{z(z-\psi)}\Big\rangle ^{Z_n}_{0,1,d} \ \right),
     \nonumber
\end{equation}
where $Q \coloneqq \mr{e}^{\tau_1}$.

The small $J$-function completely determines the restriction $\mathrm{S}(Q, z)$ of the fundamental solution to the small quantum cohomology locus. For later use, we record that using the fundamental class and divisor axiom in equation \eqref{big S-funct} one has
\begin{equation}
    \mathrm{S}(\mathbf{1}) = \frac{1}{z}\,\mathrm{J} \qquad \text{and}\qquad \mathrm{S}(H) = Q\partial_Q \mathrm{J} 
    \label{relation S(H)and S(1)}
\end{equation}
while for higher degree classes $H^a$ with $a \in \{1,\dots,n-2\}$, it holds that
\begin{equation}
    \mathrm{S}(H^{a}) = \frac{1}{\mathrm{C}^{a}_{1,a-1}}\Bigg( \partial_{\tau_1}\mathrm{S}(H^{a-1}) - \sum_{\substack{b=0 \\ b\neq a}}^{n-1} \mathrm{C}^b_{1,a-1} \, \mathrm{S}\big(H^b\big)\Bigg).
    \label{recursionS(H^a)}
\end{equation}

\section{Quantum cohomology in the \nRootSpec}
\label{sec: QuCoh}
The purpose of this section is to constrain the $Q$-dependence of the fundamental solution matrix $\mathrm{S}$. To do so, we first constrain the dependence of the small $\mathrm{J}$-function on $q$ by establishing a mirror theorem for $Z_n$. For this, we will first discuss the $I$-function of this variety. Throughout the whole section we work in the \nRootSpec\ of the equivariant parameters.

\subsection{\texorpdfstring{$I$}{I}-function}
We apply a result of Coates \cite[Thm.~3.4, 3.5]{Coates_2009} by which the $I$-function of the vector bundle $Z_n \rightarrow \mathbb{P}^{n-1}$ can be reconstructed from the one of the base
\cite[Thm.~9.5]{giventalequivariantgromovwitten}
\[
    I_{\mathbb{P}^{n-1}}(q,z) = zq^{H/z}
    \sum_{d=0}^{\infty}  \frac{1}{\prod_{i=0}^{n-1} \prod_{k=1}^{d} 
    (H + \eee_i + kz)}q^d,
\]
via a modification factor:
\begin{equation}
    I(q,z) = zq^{H/z}
    \sum_{d=0}^{\infty}  \frac{\prod_{i=0}^{n-1} \prod_{k=0}^{d-1} 
    (-H + \ddd_i - kz)}{\prod_{i=0}^{n-1} \prod_{k=1}^{d} 
    (H + \eee_i + kz)}q^d.
    \label{I function style 2}
\end{equation}

We define the logarithmic derivative $\mathrm{D}=q\frac{\partial}{\partial q}$. The $I$-function satisfies the following Picard--Fuchs equation:
\begin{equation}
    \left[ \prod_{i=0}^{n-1}(z\mathrm{D} + \eee_i) -(-1)^nq\prod_{i=0}^{n-1}(z\mathrm{D} -\ddd_i) \right]I(q,z) = 0.
    \nonumber
\end{equation}
In the \nRootSpec\ \eqref{specialization} the Picard--Fuchs equation simplifies to
\begin{equation}
    \label{eq:PFE}
    \left[ (1+(-1)^{n-1}q)(z\mathrm{D})^n +(-1)^{n-1}(\ee^n -q\dd^n)  \right]I=0.
\end{equation}
When working with the above equation we usually omit the pre-factor $zq^{H/z}$ of the $I$-function and write $I=zq^{H/z}\tilde{I}$. Then introducing the derivative operator $M \coloneqq H + z\mathrm{D}= H+zq\frac{\partial}{\partial q}$, the specialised Picard--Fuchs equation becomes
\[
    \left[ (1+(-1)^{n-1}q)M^n +(-1)^{n-1}(\ee^n -q\dd^n)  \right]\tilde{I}=0.
\]

\subsection{Mirror symmetry}
\label{sec_mirror}
The quantum Riemann--Roch theorem \cite{QuantumRiemannRoch} allows us to express the small $J$-function in terms of the $I$-function. Recall that the two must agree up to a change of coordinates $q \mapsto Q$, the mirror map, and that they both lie in the Lagrangian cone associated to the quantum cohomology of $Z_n$. Moreover, the small $J$-function is the unique element of the Lagrangian cone with the $z\rightarrow \infty$ asymptotics
\[
\mathrm{J}(Q,z) = z + \log(Q)H + \mathcal{O}(z^{-1}).
\]
See \cite[Sec.~3]{Coates_2009} for more background information. Now starting from formula \eqref{I function style 2} for the $I$-function, direct calculation shows that
\[
    I(q,z) = z + \log(q)H  +\mathcal{O}(z^{-1})
\]
meaning that the mirror map is trivial, i.e.~$Q=q$.
We have thus established the following mirror theorem.
\begin{prop}
\label{prop:mirrorsym}
The $I$-function and the small $J$-function of $Z_n$
coincide:
\[
I(q,z) = \mathrm{J}(q,z).
\]
\end{prop}

\subsection{Asymptotic expansions of the small $J$-function}
Our next goal is to constrain the $q$-dependence of the small $J$-function. For this, we write
\[
\mathrm{J}_i(q,z) = \mathrm{J}(q,z)\big|_{H=-\eee_i}
\]
for its restriction to the $T$-fixed points $P_i$ of $Z_n$, $i \in \{0,\dots,n-1\}$. Each such component admits a formal asymptotic expansion as $z \to 0$,
\begin{equation}
    \label{expansion S}
    \mathrm{J}_i \eqqcolon z \mr{e}^{u_i/z}\sum_{k\geq0}\mathrm{S}_{0i}^{[k]}z^k,
\end{equation}
since the $I$-function admits such an expansion. Here, $u_i + \eee_i\log q$ is the $i$th canonical coordinate restricted to the small quantum locus.

By \cref{prop:mirrorsym} this asymptotic expansion must satisfy the Picard--Fuchs equation \eqref{eq:PFE}. At leading order in $z$ this constraint determines the exponents $u_i$ while higher orders recursively determine the coefficients \smash{$\mathrm{S}^{\mathrlap{[k]}}{\mathstrut}_{\mathrlap{0i}\phantom{[k.}}$} up to integration constants. A priori, these are elements in $\mathbb{C}(\ee,\dd)\llbracket q \rrbracket$ but we will see that the recursive system confines them to a smaller subring. We will work with the following two rings
\begin{align*}
    R_0 & \coloneqq \
    \begin{cases}
        \mathbb{C}(\ee,\dd)\big[L^{\pm 1}\big], & \text{if $n$ is odd},\\[4pt]
        \mathbb{C}(\ee,\dd)\big[L^{\pm 1/2}\big], & \text{if $n$ is even},
    \end{cases} \qquad\qquad R \coloneqq R_0[\log L]
\end{align*}
with $L$ specified in the next \namecref{lemmwwithlogterm}. Furthermore, we write $\deg_L \mathrm{S}$ for the maximum degree in $L$ of an element $\mathrm{S}\in R_0$ or $R$ where we treat $\log L$ as a formal variable of degree zero.

\begin{lem} \label{lemmwwithlogterm}
    Let $L$ be the solution to the equation
    \begin{equation}
        \label{eq: def eq L}
        (1+(-1)^{n-1}q)L^n +(-1)^{n-1}(\ee^n -q\dd^n) = 0
    \end{equation}
    satisfying the initial condition $L\big|_{q=0}=-\ee$ and let $\zeta = \mr{e}^{2\pi \mr{i} / n}$. Then, for all $0 \le i < n$, we have
    \[ -\eee_{i} + \mr{D}u_{i} = \zeta^i\, L. \]
    Furthermore, for all $k \ge 0$
     \begin{equation}
        \label{S in R}
        \mathrm{S}_{0i}^{[k]} \in R
    \end{equation}
    and $$\deg_L \mathrm{S}_{0i}^{[k]} \leq \frac{(n-1)(2k-1)}{2}.$$
\end{lem}

\begin{rmk}
    \label{rmk: obstruction generalisation}
    This lemma is what technically obstructs us from generalising the approach pursued in \Cref{sec: QuCoh,sec: GivTel} to prove \Cref{conj: F GW via Li} beyond the \nRootSpec\ of the equivariant parameters. The specialisation is used crucially in our proof and it is not clear to the authors how it generalises. However, we expect that once an analogue of \Cref{lemmwwithlogterm} is available for general weights, all further arguments carry over almost verbatim.
\end{rmk}

\begin{proof}[Proof of \Cref{lemmwwithlogterm}]
    Instead of working with the asymptotic expansion \eqref{expansion S} we will first confine the coefficients in the expansion
    \begin{equation}
        \mathrm{J}_i \eqqcolon z \mr{e}^{u_{i}/z} \exp \left(\sum_{k\geq0}a_{0i}^{[k]}z^k \right),
        \label{expansion exp}                          
    \end{equation}
    to $R$ and only later deduce the claim \eqref{S in R}. We recall the Picard--Fuchs equation \eqref{eq:PFE}:
    \begin{equation}
        \left[(1+(-1)^{n-1}q)M^n +(-1)^{n-1}(\ee^n -q\dd^n)\right] \mathrm{J}_{i} = 0.
        \label{PicardFUchsExpansionJ}
    \end{equation}
    First, we consider how $M$ acts on the asymptotic expansion given in \cref{expansion exp}. We find 
    \[
        M \, \mathrm{J}\big|_{H=-\eee_{i}} = (L_{i} + z\mathrm{D}A)\mathrm{J}_{i}
        \nonumber
    \]
    where we define $L_{i}  = -\eee_{i}+\mathrm{D}u_{i}$ and write $A$ for  $\sum_k a_{0i}^{[k]}z^k$. This yields the following recursive formula
    \begin{equation}
        (M^n \mathrm{J})\big|_{H=-\eee_i} = P_{n,i}\, \mathrm{J}_{i},
        \nonumber
    \end{equation}
    where the sequence $\{P_{k,i}\}_{k\geq0}$ is defined recursively by
    \begin{equation}
        P_{0,i} = 1 \quad \text{and} \quad P_{k+1,i} = z\mathrm{D}P_{k,i} +(L_{i}+z\mathrm{D}A)P_{k,i}.
        \label{recursiveP}
    \end{equation}
    With this notation \cref{PicardFUchsExpansionJ} is then
    \begin{equation}
        \left[ (1+(-1)^{n-1}q)P_{n,i} +(-1)^{n-1}(\ee^n -q\dd^n)  \right]\mathrm{J}_{i}=0.
        \label{PicardFuchsRecursive}
    \end{equation}
   
     The coefficient of $z^0$ in \cref{PicardFuchsRecursive} depends only on the constant term of $P_{n,i}$, and thus we find
    \begin{equation}
        (1+(-1)^{n-1}q)L_{i}^n +(-1)^{n-1}(\ee^n -q\dd^n) =0.
        \nonumber
    \end{equation}
    Together with \(L_i|_{q=0}=-\zeta^i \ee\), this identifies \(L_i\) with the corresponding branch of \(L\). By the \(\mathbb{Z}_n\)-symmetry of the specialized equivariant weights, $ L_i = \zeta^i L$.
    
    Conversely, we can express $q$ as a rational function of $L^n$:
    \begin{equation}
        q = \frac{(-1)^n L^n - \ee^n}{L^n - \dd^n}.
        \nonumber
    \end{equation}
    Differentiating the expression for $q$ with respect to $L = L_0$ yields
    \begin{equation}
        \mathrm{D}L = q\frac{\partial L}{\partial q}=\frac{((-1)^nL^n-\ee^n)(L^n-\dd^n)}{n L^{n-1}(\ee^n-(-1)^n\dd^n)},
        \nonumber
    \end{equation}
    which in particular lies in $R$.
    
    Subtracting the defining equation of $L_i$ from \cref{PicardFuchsRecursive} we find
    \[
    (1+(-1)^{n-1}q)(P_{n,i}-L_i^n) = 0,
    \]
    and thus $P_{n,i} = L_i^n$.
    Expanding in powers of $z$, the coefficient of $z^{k+1}$ has the form
    \begin{equation}
        \label{recursion_coeff pre}
        nL_i^{n-1} \mathrm{D} a_{0i}^{[k]} + E_{k,i} =0,
    \end{equation}
    where $E_{k,i}$ depends on $a_{0i}^{[0]}, \dots, a_{0i}^{[k-1]}, L_i$ and their $\mathrm{D}$-derivatives. To produce $a_{0i}^{[k]} $ at this order, one has to choose the factor $z\mathrm{D}A$ once in the $n$-fold recursion, and all other factors must be $L_i$. Thus the coefficient is $nL_i^{n-1}$. Furthermore, every term in the coefficient contains at least one factor of $\mathrm{D}$, since one must choose either $z\mathrm{D}A$ or $z\mathrm{D}$ to produce a positive power of $z$. Since $\mathrm{D} f(L_i) = df/d{L_i}\cdot  \mathrm{D}L_i$ we can factor out the global factor $\mathrm{D}L_i$ to obtain
    \begin{equation}
        \frac{d}{dL_i} a_{0i}^{[k]} = -\frac{1}{nL_i^{n-1}} \tilde{E}_{k,i}.
        \label{recursion_coeff}
    \end{equation}
    For $k=0$ we find
    \[
        \frac{d}{dL_i}a_{0i}^{[0]} = -\frac{n-1}{2}\frac{1}{L_i},
    \]
    which has general solution
    \begin{equation}
        \label{a0}
        a^{[0]}_{0i} = \frac{n-1}{2}\log \left(\frac{c_{1,i}}{L_i} \right),
    \end{equation}
    for some constant $c_{1,i}$. Thus
    \begin{equation}
        \label{S0}
        S_{0i}^{[0]} = \exp\big( a_{0i}^{[0]} \big)
    \end{equation}
    lies in $R$. Using that $\tilde{E}_{k,i}\in R$ and that $R$ is closed under integration, it recursively follows from equation \cref{recursion_coeff} that every $a_{0i}^{[k]}\in R$.

    Let us now use the recursion to prove that
    \begin{equation}
        \label{deg a k}
        \deg_L a_{0i}^{[k]} \leq k(n-1).
    \end{equation}
    For this, it will be convenient to extend $\deg_L$ to the variable $z$ by formally assigning it the degree $\deg_L z = 1-n$. This means we want to show that \smash{$\deg_L a^{\mathrlap{[k]}}{\mathstrut}_{0i}~ \, z^k \leq 0$} for all $k \geq 0$. We will prove this by induction on $k$. Note that by \eqref{a0} we already know that the claim holds for $k=0$. Recall that we formally assigned $\deg_L \log L =0$. So suppose the degree bound holds for all $0\leq j <k$. Then since in the recursion \eqref{recursiveP} the operator $z \mathrm{D}$ may at most raise the degree by one and $z^{k+1} E_{k,i}$ on the right-hand side of \eqref{recursion_coeff pre} only receives contributions from the truncation of $A$ at $z^{k-1}$, which by assumption has degree at most $0$, we have $\deg_L z^{k+1} E_{k,i} \leq n$. Thus, by our recursion \eqref{recursion_coeff pre}, we must have $\deg_L \mathrm{D} a^{\mathrlap{[k]}}{\mathstrut}_{0i}~ \, \leq k(n-1) + n$ from which the claimed maximum degree \eqref{deg a k} for the $k$th coefficient follows. 
    
    Comparing the expansions \eqref{expansion S} and \eqref{expansion exp}, the coefficients \smash{$\mathrm{S}^{\mathrlap{[k]}}{\mathstrut}_{0i}~$} are related to the coefficients \smash{$a^{\mathrlap{[k]}}{\mathstrut}_{0i}~$} by
    \begin{equation}
        \label{Bell expansion}
        \mathrm{S}_{0i}^{[k]} = \mathrm{S}_{0i}^{[0]} \, \frac{B_k\big(1!\,a_{0i}^{[1]},2!\,a_{0i}^{[2]},\dots,k!\,a_{0i}^{[k]}\big)}{k!},
    \end{equation}
    where $B_k$ is the complete Bell polynomial of order $k$.
    We conclude that
    \begin{equation*}
        \mathrm{S}_{0i}^{[k]} \in R.
    \end{equation*}
    Regarding the $L$-degree bound, it follows from equations \eqref{a0} and \eqref{S0} that $\deg_L \mathrm{S}^{\mathrlap{[0]}}{\mathstrut}_{0i}~ = -(n-1)/2$. Moreover, the Bell polynomial factor on the right-hand side of equation \eqref{Bell expansion} has $L$-degree at most $k(n-1)$ since each term in $$\exp \sum_{k>0} a_{0i}^{[k]} z^k$$ has $L$-degree at most $0$ by \eqref{deg a k}. Hence, combining the degrees of both factors on the right-hand side of equation \eqref{Bell expansion}, we indeed arrive at the claimed bound
    \begin{equation*}
        \deg_L \mathrm{S}_{0i}^{[k]} \leq \frac{(n-1)(2k-1)}{2}. \qedhere
    \end{equation*}
\end{proof}


\subsection{A strengthening via Givental's equivariant mirror}
\label{sec:GivMirror}
To strengthen our last result and confine the asymptotic expansion of the $J$-function to $R_0\subset R$ we require its presentation as an asymptotic expansion of an oscillatory integral arising from Givental's equivariant mirror \cite{giventalequivariantgromovwitten}. An introduction to the methods used can be found in \cite[Appendix]{CrepantresolutionC3/Z3}.

Givental introduced the equivariant mirror for $Z_n$ as the Landau--Ginzburg potential
\begin{equation}
    W(w_0,\dots,w_{2n-1}) = \sum_{i=0}^{2n-1}w_i - \sum_{i=0}^{n-1} \eee_i \log(w_i) -\sum_{i=n}^{2n-1}\ddd_{i-n} \log(w_i),
    \nonumber
\end{equation}
defined on the family of affine varieties
\begin{equation}
    M_{q} = \left\{(w_0,\dots,w_{2n-1})\in \mathbb{C}^{2n}:\prod_{i=0}^{n-1} w_i = q \prod_{i=n}^{2n-1}w_i\right\}.
    \label{nonlinearconstraint}
\end{equation}
We have the associated oscillatory integral
\begin{equation}
    \mathcal{I}(z,\Gamma)= \int_{ \Gamma \subset M_{q}} \mr{e}^{W/z} \omega,
    \label{osciIntegral}
\end{equation}
where $\Gamma \subset M_{q}$ is a suitable $(2n-1)$-cycle and $\omega$ denotes the meromorphic volume form:
\[
\omega = \frac{d\log(w_0) \wedge \ldots \wedge d\log(w_{2n-1}) }{d \log (q)}.
\]

Let $\text{Asym}_{\text{cr}_i}(\mr{e}^{W/z}\omega)$ denote the formal asymptotic expansion of the oscillatory integral \cref{osciIntegral} at a critical point $\text{cr}_i$ neglecting exponentially sub-leading terms. By \cite[Proposition 6.9]{OscillatoryIntegralRelation}, there is a relation between the $I$-function and the oscillatory integral. In our setting, it reads
 \begin{equation}
    \mr{e}^{W(\text{cr}_i)/z}\cdot\text{Asym}_{\text{cr}_i}\left(\mr{e}^{W/z}\omega\right) =
    \mr{e}^{u_i/z}\cdot \left( 1+\frac{\mathrm{S}_{0i}^{[1]}}{\mathrm{S}_{0i}^{[0]}}z+\frac{\mathrm{S}_{0i}^{[2]}}{\mathrm{S}_{0i}^{[0]}}z^2+\dots \right).
    \label{relation oscillatory Integraland expansion}
\end{equation}
This presentation allows us to prove the following strengthening of \Cref{lemmwwithlogterm}.

\begin{lem} \label{lem:no log}
    For $k \geq 0$ and $i \in \{0,\dots,n-1 \}$ it holds that
    \begin{equation}
        \mathrm{S}_{0i}^{[k]} \in
        R_0
    \end{equation}
\end{lem}

\begin{proof}
We adapt an argument of Lho used in \cite[Lemma 7]{LhoLocalP1}. The saddle point method allows one to compute the asymptotic behaviour of the oscillatory integral \eqref{osciIntegral}. We let $\operatorname{H}_{i} \coloneqq\operatorname{Hess}W(\text{cr}_i)$ denote the Hessian at a critical point. The formal asymptotic expansion of $\text{Asym}_{\text{cr}_i}(\mr{e}^{W/z}\omega)$ then takes the explicit form:
\begin{equation}
    \text{Asym}_{\text{cr}_i}(\mr{e}^{W/z}\omega) =\frac{1}{\sqrt{\det(\operatorname{H}_{i})}} \left.
    \exp\!\left(-\frac{z}{2}\,\nabla^{T} \operatorname{H}_{i}^{-1} \nabla       \right) \mr{e}^{W_{\ge 3}(w)/z} \right |_{w=\mathrm{cr}_i},
    \label{expansionOscillatoryIntegral}
\end{equation}
where $\nabla$ denotes the gradient operator with respect to the variables $w$ and the exponential operator is defined by its formal power series, acting on the terms to the right of it.
Here $W_{\geq 3}$ denotes the Taylor remainder of $W$ at $\text{cr}_i$, given by discarding all terms of degree two and less:
\[
W_{\geq 3}(w)=W(w) - W(\text{cr}_i) - \frac{1}{2} (w-\text{cr}_i)^{\mr{t}} \operatorname{H}_{i}(w-\text{cr}_i).
\]

To determine the critical points of $W$ on $M_{q}$, we apply the method of Lagrange multipliers. We rewrite the condition \eqref{nonlinearconstraint} as
 \begin{equation}
     F \coloneqq \log\left(\; \prod_{i=0}^{n-1}w_i\right) - \log\left( \;\prod_{i=n}^{2n-1}w_i\right) - \log(q) = 0.
     \nonumber
 \end{equation}
Denoting the Lagrange multiplier by $\mu$, the criticality conditions read
\begin{equation}
    \frac{\partial W}{\partial w_i} = \mu \frac{\partial F}{\partial w_i} .
    \nonumber
\end{equation}
Evaluating the differentials gives
\begin{equation}
    w_i = \begin{cases}
    \mu + \zeta^i \ee, & i \in \{0, \dots, n-1\}, \\
    -\mu + \zeta^i \dd, & i \in \{n, \dots, 2n-1\}.
    \end{cases}
\label{134}
\end{equation}
Substituting these into condition \eqref{nonlinearconstraint}, we obtain
\begin{equation}
    \prod_{i=0}^{n-1}(\mu+\zeta^i \ee) = q\prod_{i=n}^{2n-1}(-\mu+\zeta^i \dd).
    \nonumber
\end{equation}
By the symmetric properties of the $n$th roots of unity, this simplifies to 
\begin{equation}
    \mu^n+(-1)^{n-1}\ee^n=q(-1)^n(\mu^n -\dd^n).     
    \nonumber
\end{equation}
This is precisely the defining equation of $L$, and each solution for $\mu$ corresponds to a critical point. In what follows, we choose the branch with $L|_{q=0}=-\ee$ and denote the corresponding critical point by $\text{cr}_0$. Setting $\mu = L$ in expression \eqref{134} yields
\[
\begin{aligned}
w_i &= L + \zeta^i \ee, && \text{for } i \in \{0,\dots,n-1\}, \\
w_i &= -L + \zeta^i \dd, && \text{for } i \in \{n,\dots,2n-1\}.
\end{aligned}
\]

The determinant of the Hessian of $W$ with respect to the logarithmic coordinates $(\log( w_0), \allowbreak \dots, \allowbreak \log( w_{2n-1}))$ at the critical point $\text{cr}_0$ is then
\begin{equation}
    \det\operatorname{Hess}W(\text{cr}_0)  =
    nL^{n-1}(\ee^n+ (-1)^{n-1}\dd^n).
    \label{detHess}
\end{equation}

To ensure that the formal asymptotic expansion \eqref{expansionOscillatoryIntegral} lies in the claimed ring,
we must verify that the edge weights (arising from $\mr{H}_i^{-1}$) and 
vertex weights (arising from $W_{\geq 3}$) introduce no terms $\log(L)$.

In the logarithmic coordinates $x_i = \log w_i$, the potential takes the form
\begin{equation}
    W = \sum \mr{e}^{x_i} - \sum \eee_i x_i - \sum \ddd_i x_i.
    \label{potentiallogcoords}
\end{equation}
Because the logarithmic terms in equation \eqref{potentiallogcoords} are linear in $x_i$, 
they are annihilated by all derivatives of order two or higher. 
Thus, the Hessian $H_i$ and the Taylor remainder $W_{\geq 3}$ 
depend only on derivatives of the exponential part:
\[
    \frac{\partial^k}{\partial x_i^k} \mr{e}^{x_i} = \mr{e}^{x_i} = w_i.
    \nonumber
\]

When evaluated at the critical point $\mathrm{cr}_0$, 
these derivatives simply yield the coordinate values $w_i$, 
which are linear in $L$, $\ee$, and $\dd$. 
Thus, all vertex weights and the entries of $H_i$ 
lie in the ring $\mathbb{C}(\ee,\dd)[L]$.

Inverting the Hessian amounts to dividing its adjugate matrix by the determinant \eqref{detHess}, introducing 
only a factor $L^{-(n-1)}$. This completes the verification that \eqref{expansionOscillatoryIntegral}
lies in $R_0$.

By relation \eqref{relation oscillatory Integraland expansion} it follows that the coefficients $\mathrm{S}_{0i}^{[k]}/\mathrm{S}_{0i}^{[0]}$ also lie in this ring. This proves the claim for both parities of $n$.
\end{proof}

\subsection{Identification of all structure constants}
Remember that the goal of this section is to constrain the $q$-dependence of the fundamental solution matrix. So far we have only proved this for its first column. To prepare for extending our result on the first column to the remaining ones, we first determine the structure constants of the small quantum product.

\begin{lem}
    The structure constants of the small quantum product of $Z_n$ satisfy
    \[
    \mathrm{C}^c_{ab} = \delta_{c,a+b},
    \]
    when $a+b<n.$
    \label{lem:structure const}
\end{lem}

\begin{proof}
The small quantum product is given by
\[
H^a \circ_{\tau} H^b = \sum_{c,c'} \sum_{d } \eta^{cc'}\,\langle H^a, H^b, H^c \rangle_{0,3,\beta}\,  H^{c'}\ q^{d}.
\]
Hence, we need to evaluate the correlator
\begin{equation}
    \label{3pt correlator}
    \langle H^a, H^b, H^{n-1-c}  \rangle_{0,3,d}^{Z_n} = \int_{[\Mbar_{0,3}(Z_n,d)]_{T}^{\mathrm{vir}}}\mathrm{ev}_1^*(H^a)\mathrm{ev}_2^*(H^b)\mathrm{ev}_3^*(H^{n-1-c}).
\end{equation}
By a dimension count, the integral takes values in the $(a+b-c-n)$-graded piece of localised equivariant Chow cohomology of a point. In other words, the integral is a rational function in $\ee$ and $\dd$ with homogeneous numerator and denominator whose degree difference is $a+b-c-n$. Since by \Cref{lem: s t dependence} the dependence on these variables must be through $\dd^n$ and $\ee^n$, the correlator vanishes unless
\[
    a+b-c \equiv 0 \mod n .
\]
Since $0\leq c\leq n-1$ and we assume that $a+b < n$, we get
\[
    -n+1\leq a+b-c\leq n-1
\]
The only $c$ satisfying both is
\[
    a+b=c.
\]
However, in this case the degree of the correlator \eqref{3pt correlator} in $\ee$ and $\dd$ must be $-n$, ie.~it is strictly negative. Since by the properness of the moduli space in positive degree $d>0$ the dependency on these variables must be polynomial (\Cref{lem: s t dependence}), we conclude that
\begin{equation*}
    \langle H^a, H^b, H^{n-1-c}  \rangle_{0,3,d}^{Z_n} = 0
\end{equation*}
whenever $a+b < n$ and $d>0$.

The \namecref{lem:structure const} then follows from the evaluation of the classical three-point correlator
\begin{equation*}
\mathrm{C}^{a+b}_{ab} = \eta^{a+b,n-1-a-b}\langle H^a, H^b, H^{n-1-a-b}\rangle = \eta^{a+b,n-1-a-b}\eta_{a+b,n-1-a-b} = 1. \qedhere
\end{equation*}
\end{proof}


\subsection{Canonical to flat frame}
As a further preparation for constraining the $q$-depen\-dence of the fundamental solution matrix, we determine the change of basis matrix from the canonical to the flat frame in this subsection. Recall that the flat frame is given by the basis $\{H^a\}_{a\in\{0,\dots,n-1\}}$ while the canonical frame is given by the idempotents $\{\mathbf{e}_i\}_{i\in\{0,\dots,n-1\}}$ of the quantum product. We define the transition matrix 
\[
    \Phi_{ai} \coloneqq (\partial_a u_i)\big|_{\tau_b=\delta_{b,1}\log(q)}
\]
from the canonical frame to the flat frame, restricted to the small quantum cohomology locus. The normalisation factor is defined by
\[
    \Delta_i^{-1} \coloneqq \eta(\mathbf{e}_i, \mathbf{e}_i),
\]
and we write
\[
    \Psi_{ai} \coloneqq \Phi_{ai} \Delta_i^{-1/2}
\]
for the transition matrix from the normalised canonical frame to the flat frame.

\begin{prop}
    \label{prop:PhiPsiDelta}
    We have
    \begin{equation*}
        \Phi_{ai} = L_i^a, \qquad \Delta_i = n L_i^{\,n-1} \bigl( \ee^n - (-1)^{n}\dd^n \bigr).
    \end{equation*}
    In particular, $\Phi_{ai}, \Delta_i \in 
    R_0$.
\end{prop}

\begin{proof}
    By idempotency we have
    \[
        \Phi_{ai}\Phi_{bi} = \sum_{c=0}^{n-1}\mathrm{C}^c_{ab}\Phi_{ci},
    \]
    From the string axiom and \cref{int_prod} we get $\Phi_{0i} = 1$ and that $\Phi_{1i} = L_i$. Inserting the structure constants imposes $\Phi_{(a+1) i}= \Phi_{ai}\Phi_{1i}$.
    We conclude that
    \[
        \Phi_{ai} = L_i^a,
    \]
    for $a \in \{0,\dots,n-1\}$. This allows us to determine the normalisation term $\Delta_i^{-1}$ of the idempotent basis:
    \begin{align*}
        \Delta_i^{-1} 
        &= \eta(\mathbf{e}_i, \mathbf{e}_i) = \eta(\mathbf{e}_i, \mathbf{1}) 
        = \sum_{c=0}^{n-1} (\Phi^{-1})_{ic} \, \eta(H^c, \mathbf{1}) \\
        &= (\Phi^{-1})_{i,n-1} \, \eta_{n-1,0} 
        = \frac{1}{n L_i^{\,n-1} \bigl( \ee^n - (-1)^{n}\dd^n \bigr)}.\qedhere
    \end{align*}
\end{proof}

\subsection{Finite generation of the fundamental solution}
We are now equipped to confine the coefficients of the $z \rightarrow 0$ asymptotic expansion of the fundamental solution matrix to $R_0$. For this, note that equation \eqref{recursionS(H^a)} expresses $\mathrm{S}(H^a)$ for $a \in \{0,\dots,n-1 \}$ in terms of the small $J$-function and the structure constants. Substituting the structure constants found in \cref{lem:structure const} into \cref{recursionS(H^a)} yields
\begin{equation}
    \mathrm{S}(1) = \frac{1}{z}\mathrm{J}, \quad
    \mathrm{S}(H) = \mathrm{D}\mathrm{J}, \quad
    \cdots \quad
    \mathrm{S}(H^a) = z^{a-1}\mathrm{D}^a\mathrm{J}.
    \label{relations all}
\end{equation}
Since the small $J$-function allows an asymptotic expansion when restricted to the fixed points as $z \to0$, we obtain similar expansions for the fundamental solutions $\mathrm{S}(H^a)$. We consider the expansions
\[
\mathrm{S}_{ai} \coloneqq \sum_{k=0}^{\infty} \mathrm{S}^{[k]}_{ai}\, z^{k} \coloneqq \mr{e}^{-u_i/z} \, \mathrm{S}(H^a)\big|_{H=-\eee_i}.
\]

\begin{prop} \label{prop:higher S expansion in ring}
    For all $a,i\in \{  0,1,\dots,n-1\}$ and $k\geq 0$ the coefficients $\mathrm{S}^{[k]}_{ai}$ lie in $R_0$. Furthermore,
    \begin{equation*}
        \deg_L \mathrm{S}^{[k]}_{ai} \leq a+ \frac{(n-1)(2k-1)}{2}.
    \end{equation*}
\end{prop}
\begin{proof}
    We already proved the claim for $a=0$ in \Cref{lemmwwithlogterm} and \ref{lem:no log}. The case $a>0$ is reduced to this base case via the tower \eqref{relations all} observing that the operator $\mr{D}$ preserves $R_0$ and $z\mr{D}$ raises the $L$-degree by at most one if, as before, we formally assign $\deg_L z = 1-n$.
\end{proof}

\section{Higher-genus reconstruction in the \nRootSpec}
\label{sec: GivTel}
In this section we will use the Givental--Teleman reconstruction theorem \cite{giventalconjecture,Telemanproof} to transport our earlier results on the fundamental solution matrix to the higher genus Gromov--Witten theory of $Z_n$ in order to prove \Cref{thm: F GW via Li}\ref{item: thm B i}. As before, throughout the whole section we will work in the \nRootSpec\ of the equivariant parameters.

\subsection{Identification of the $R$-matrix}
The central algebraic object in the reconstruction procedure of Givental--Teleman is the $R$-matrix. To compute this matrix explicitly, we characterize it as the particular fundamental solution $R_j(\tau,z)$ to the Dubrovin connection satisfying:
\begin{description}[labelindent=1em]
    \item[Asymptotics] We have $R_j(\tau,z)(H^a) = \mr{e}^{u_j/z}\,R_{aj}(\tau,z)$ with $R_{aj} \in \mathbb{C}(\epsilon,\delta) [\![z]\!]$ for fixed $\tau$.
    
    \item[Symplecticity] We have $\sum_{a,b} R_{ai}(\tau,z)\,\eta^{ab}\,R_{bj}(\tau,-z) = \delta_{i,j}$.
    
    \item[Normalisation] The small quantum cohomology restriction $R(\tau)\mapsto \mathrm{R}(q)$ satisfies
    \[
    \mathrm{R}_{aj}\big|_{q=0} = \Psi_{aj}\big|_{q=0}\, \mr{e}^{-c_j(z)},
    \] 
    where 
    \begin{equation}
        c_j(z) \coloneqq 
        \sum_{k=1}^{\infty}
        \frac{B_{2k}}
        {2k(2k-1)} \, z^{2k-1} \sum_{i=0}^{2k-2} \big(w_j^{i}\big)^{1-2k}.
        \label{cj(z)}
    \end{equation}
    Here, $(w_j^{i})_{1\leq i \leq 2n}$ denote the tangent $T$-weights at the fixed point $P_j \in Z_n$.
\end{description}
Such a solution exists since the quantum cohomology of $Z_n$ is semi-simple. We start by constructing the $R$-matrix of $Z_n$ from the fundamental solution $S(\tau,z)$ to the Dubrovin connection introduced in \Cref{sec: J fct}. It satisfies the symplecticity condition and we showed that it admits the required asymptotics. Since, moreover, $\mathrm{f}_j^*\mathrm{S}(q,z)(H^a) \rightarrow \mathrm{f}_j^*H^a |_{q=0} = \Psi_{aj}|_{q=0}$ as $q \rightarrow 0$, we deduce that
\begin{equation}
    \mathrm{R}_{aj} = \mathrm{S}_{aj} \, \mr{e}^{-c_j(z)},
    \label{R-S matrix}
\end{equation}
where $c_j(z)$ is defined as in \cref{cj(z)}. The following proposition is then an immediate consequence of \Cref{prop:higher S expansion in ring}.

\begin{prop}
    \label{prop:R in ring}
    All coefficients in the expansion $\mathrm{R}_{aj} = \sum_{k \geq 0} \mathrm{R}_{aj}^{[k]} \, z^k$ lie in $R_0$ and
    \begin{flalign*}
        && \deg_L \mathrm{R}^{[k]}_{aj} \leq a+\frac{(n-1)(2k-1)}{2}. && \qed
    \end{flalign*}
\end{prop}

Later we will also require the following symmetry of the $R$-matrix.

\begin{prop}
    \label{lem: R q sym}
    Under the substitution
    \begin{equation}
        \label{eq: q q inv subs}
        q \longmapsto q^{-1}, \qquad (\eee_0,\ldots,\eee_{n-1}) \longleftrightarrow (\ddd_0,\ldots,\ddd_{n-1})
    \end{equation}
    the entries of the $R$-matrix transform as $\mathrm{R}_{aj} \mapsto (-1)^a \mathrm{R}_{aj}$.
\end{prop}

\begin{proof}
    It suffices to prove the symmetry for the entries $\mathrm{R}_{0j}$ since all others are recovered via \eqref{relations all}. To prove the symmetry for the identity column, we once more require Givental's equivariant mirror. Recall from \Cref{sec:GivMirror} that the entries $\mathrm{S}_{0j}$ are proportional to a formal asymptotic expansion
    \begin{equation}
        \label{eq: Asym LG}
        \text{Asym}_{\text{cr}_j}\left(\mr{e}^{W/z}\omega\right)
    \end{equation}
    of an oscillatory integral over a cycle in a family $M_q$ of affine varieties with a suitable Landau--Ginzburg potential $W$ and volume form $\omega$. To prove our proposition we only require two facts about this asymptotic expansion: First, by Stirling's approximation of the Gamma function we have
    \begin{equation*}
        \text{Asym}_{\text{cr}_j}\left(\mr{e}^{W/z}\omega\right) \Big|_{q=0} = \Delta_j^{-1/2}\big|_{q=0} \,\re^{-c_j(z)}.
    \end{equation*}
    Together with \eqref{relation oscillatory Integraland expansion} this implies that
    \begin{equation}
        \label{eq: R to asymp}
        \mathrm{R}_{0j}=\text{Asym}_{\text{cr}_j}\left(\mr{e}^{W/z}\omega\right).
    \end{equation}
    The second property we require about Givental's equivariant mirror is that $M_q$ and $W$ are invariant under the substitution \eqref{eq: q q inv subs} and swapping $w_i\leftrightarrow w_{i+n}$ implying that the same symmetry holds for the left-hand side of \eqref{eq: R to asymp}.
\end{proof}

\subsection{Decorated stable graphs and their weights}
The Givental--Teleman reconstruction theorem expresses a Gromov--Witten correlator of $Z_n$ as a weighted sum over certain decorated graphs.

\subsubsection{Definition}
\label{def:decorated-stable-graph}
An $\boldsymbol{n}$\textbf{-decorated stable graph}
\[
    \Gamma = (V_\Gamma,\, H_\Gamma,\, E_\Gamma,\, L_\Gamma,\, D_\Gamma,\, g,\, d,\, i,\, v,\, q,\, a)
\]
consists of the following data:
\begin{itemize}
    \item $V_\Gamma$: a finite set of vertices, equipped with
    \begin{itemize}
        \item a genus assignment $g : V_\Gamma \to \mathbb{Z}_{\geq 0}$, recording the geometric genus at each vertex;
        \item a marking assignment $d : V_\Gamma \to \{0, 1, \ldots, n-1\}$, recording the cohomological state inserted at each vertex;
    \end{itemize}
    
    \item $H_\Gamma$: a finite set of half-edges, equipped with
    \begin{itemize}
        \item an involution $i : H_\Gamma \to H_\Gamma$;
        \item an attaching map $v : H_\Gamma \to V_\Gamma$, recording the vertex to which each half-edge is attached;
    \end{itemize}
    
    \item $E_\Gamma$: the set of edges, defined as the set of $2$-orbits $e = \{h_1, h_2\}$ of the involution $i$; self-edges (loops) are permitted;
    
    \item $L_\Gamma \sqcup D_\Gamma$: a partition of the set of fixed points of $i$ (equivalently, the set of leaves) into
    \begin{itemize}
        \item $L_\Gamma$: a set of marked leaves, corresponding to the marked points;
        \item $D_\Gamma$: a set of dilaton leaves, corresponding to dilaton insertions;
    \end{itemize}
    
    \item a leg marking assignment $q : L_\Gamma \to \{0, 1, \ldots, n-1\}$, recording the cohomological state at each marked leaf;
    
    \item an ancestor assignment $a : H_\Gamma \to \mathbb{Z}_{\geq 0}$, recording the $\psi$-class degree at each half-edge.
\end{itemize}
These data are required to satisfy the following two conditions:
\begin{enumerate}
    \item The pair $(V_\Gamma, E_\Gamma)$ defines a connected graph.
    
    \item For each vertex $v \in V_\Gamma$, letting $H(v)$ denote the set of half-edges incident to $v$, the stability condition
    \[
    2g(v) - 2 + |H(v)| > 0
    \]
    is satisfied.
\end{enumerate}
The genus of $\Gamma$ is defined as
\[
g(\Gamma) \;=\; h^1(\Gamma) + \sum_{v \in V_\Gamma} g(v),
\]
where $h^1(\Gamma) = |E_\Gamma| - |V_\Gamma| + 1$ is the first Betti number of the underlying graph. We denote by \smash{$\mathcal{G}^{\mathrlap{[n]}}{\mathstrut}_{g,m}$} the set of isomorphism classes of $n$-decorated stable graphs of genus $g$ with $m$ marked legs.

\subsubsection{Graph weights}
With the graph framework in place, the next step is to specify how the geometry of the target influences the reconstruction formula. We make this precise by assigning explicit weights in terms of the canonical $R$-matrix to the vertices, edges and legs.
Given a $n$-decorated stable graph $\Gamma$ we assign the following weights to each leaf, vertex and edge. For $ i,j \in \{0,\ldots,k\} $ and $ l,m \in \mathbb{Z}_{\ge 0} $, we define the following:
\begin{itemize}
    \item
     The \emph{dilaton factor} is
    \begin{equation}
    D_{i,l}
    \coloneqq
    [z^{l-1}]\bigl(\Psi_{0i} - \mathrm{R}_{0i}(-z)\bigr),
    \label{dilaton}
    \end{equation}
    \item The \emph{marking factor} is
    \begin{equation}
    L^{[a]}_{i,l}
    \coloneqq
    [z^l] \mathrm{R}_{ai}(-z).      
    \label{leaf}
    \end{equation}
    
    \item
    The \emph{edge factor} is given by
    \begin{equation}
    E^{ij}_{l,m}
    \coloneqq
    [z^l w^m]
    \left(
    \frac{
    \delta_{i,j}
    -
    \sum_{ab}
    \mathrm{R}_{ai}(-z)\,\eta^{ab}\,\mathrm{R}_{bj}(-w)
    }{z+w}
    \right).         
    \label{edge}
    \end{equation}
    
    \item
    Each \emph{leaf} $ h \in D_{\Gamma} \sqcup L_{\Gamma} $ contributes
    \begin{equation}
    \mathrm{Cont}(h)
    \coloneqq
    \begin{cases}
    D_{q(v(h)),\,a(h)}, & h \in D_{\Gamma}, \\[6pt]
    L^{[q(h)]}_{q(v(h)),\,a(h)}, & h \in L_{\Gamma} .
    \end{cases}         
    \nonumber
    \end{equation}
    
    \item
    Each \emph{edge} $ e = \{h_1,h_2\} \in E_{\Gamma} $, is attributed
    \begin{equation}
    \mathrm{Cont}(e)
    \coloneqq
    E^{\,q(v(h_1)),\,q(v(h_2))}_{\,a(h_1),\,a(h_2)}.
    \nonumber
    \end{equation}
    
    \item
    For a fixed labelling $\{h_1,\dots,h_{|H(v)|}\}$ of the half-edges incident to a vertex $ v \in V_{\Gamma} $ we define
    \begin{equation}
    \mathrm{Cont}(v)  =   
    (\Delta_{q(v)})^{\frac{2g(v)-2+|H(v)|}{2}}       
    \int_{[\Mbar_{g(v),\,|H(v)|}]}
    \prod_{k=1}^{|H(v)|}
    \psi_k^{\,a(h_k)} .
    \nonumber
    \end{equation}
\end{itemize}

The following is a consequence of \Cref{prop:PhiPsiDelta,lem: R q sym,prop:R in ring}.

\begin{prop}
    \label{prop:dilaton leaf edge vertex}
    The dilaton factor $D_{i,l}$, the marking factor $L_{i,l}^{[a]}$ and the edge factor $E_{l,m}^{i,j}$ of $Z_n$ lie in $R_0$. The dilaton and edge factor are invariant under the simultaneous exchange of $q\rightarrow q^{-1}$ and $\eee_i \leftrightarrow \ddd_i$, $i\in\{0,\ldots,n-1\}$, while marking factors transform as $L_{i,l}^{[a]} \mapsto (-1)^a L_{i,l}^{[a]}$. Furthermore,
    \begin{gather*}
        \deg_L D_{i,l} \leq \frac{(n-1)(2l-3)}{2},\quad \deg_L L_{i,l}^{[a]} \leq a + \frac{(n-1)(2l-1)}{2},\\[0.7em]
        \deg_L E_{l,m}^{i,j} \leq (n-1)(l+m+1).
    \end{gather*}
    For every vertex $v$ we have $\mathrm{Cont}(v) \in R_0[\Delta^{\mathrlap{1/2}}{\mathstrut}_{\mathrlap{i}\phantom{1/2}}]$ and
    \begin{flalign*}
        && \deg_L \mathrm{Cont}(v) \leq \frac{(n-1)(2g(v)-2+|H(v)|)}{2}. && \qed
    \end{flalign*}
\end{prop}


\subsection{Finite generation via Givental--Teleman reconstruction}
Using the fact that the $R$-matrix is finitely generated, we extend our previous results about restricted $q$-dependence from genus-zero to arbitrary genus. Recall our notation for the fixed-genus generating series of $T$-equivariant Gromov--Witten invariants of $Z_n$ from the introduction:
\[
    F{}^{Z_n,T}_{g}(q) = \sum_{d > 0} q^d  \int_{[\Mbar_g(Z_n,d)]^{\vir}_T}1.
\]
Note that we only sum over strictly positive degrees $d>0$. Via the divisor equation, the $m$th logarithmic $q$-derivative of this series can be written in our correlator notation as
\begin{equation}
    \label{F g vs correlator}
    \mathrm{D}^m \,F^{Z_n,T}_{g} = \bigl\langle\!\!\bigl\langle H, \ldots, H \bigr\rangle\!\!\bigr\rangle^{Z_n}_{g,m} - \int_{[\Mbar_g(Z_n,0)]^{\vir}_T}1.
\end{equation}
It therefore suffices to analyse the right-hand side correlator in order to prove \Cref{thm: F GW via Li}\ref{item: thm B i}.

\begin{prop} \label{thrm:finitley generated n}
    Let $g,m\geq 0$ with $2g-2+m > 0$. Then
    \[
    \bigl\langle\!\!\bigl\langle H, \ldots, H \bigr\rangle\!\!\bigr\rangle^{Z_n}_{g,m} \in R_0
    \]
    Furthermore,
    \begin{equation}
        \label{eq: F GW L deg}
        \deg_L \bigl\langle\!\!\bigl\langle H, \ldots, H \bigr\rangle\!\!\bigr\rangle^{Z_n}_{g,m} \leq (n-1)(4g-4+m)+m
    \end{equation}
    and the correlator has parity $(-1)^m$ under the simultaneous exchange of $q\mapsto q^{-1}$ and $\eee_i \leftrightarrow \ddd_i$, $i\in\{0,\ldots,n-1\}$.
\end{prop}
\begin{proof}
    Let us more generally consider correlators with arbitrary primary insertions. For those the Givental--Teleman formula \cite{giventalconjecture,Telemanproof} states that
    \[
    \bigl\langle\!\!\bigl\langle
    H^{a_1}, \ldots, H^{a_m}
    \bigr\rangle\!\!\bigr\rangle^{Z_n}_{g,m} =            \sum_{\Gamma \in \mathcal{G}^{[n]}_{g,m}}\! \frac{1}{|\mathrm{Aut}(\Gamma)|} 
    \prod_{e \in E_\Gamma} \mathrm{Cont}(e) \!\!\!
    \prod_{h \in L_\Gamma \sqcup D_\Gamma}\!\!\! \mathrm{Cont}(h) \!
    \prod_{v \in V_\Gamma}     \mathrm{Cont}(v).
    \]
    By \cref{prop:dilaton leaf edge vertex}, every term of the right-hand side sum lies in $R_0[\Delta_0^{1/2}]$ and satisfies the anticipated symmetry. We can deduce that the overall sum lies in $R_0$ as Gromov--Witten invariants are rational in $\ee$ and $\dd$ and \smash{$\Delta^{\mathrlap{1/2}}{\mathstrut}_{\mathrlap{i}\phantom{1/2}} \propto \sqrt{ \ee^n - (-1)^{n}\dd^n}$}.
    
    Moreover, by \cref{prop:dilaton leaf edge vertex}, the $L$-degree of a term associated to a decorated stable graph $\Gamma$ is at most
    \begin{equation*}
        m+ \frac{n-1}{2} \left( \sum_v \big(2g(v)-2+|H(v)|\big)  - |L_{\Gamma}| + 2|E_{\Gamma}| - 3|D_{\Gamma}| + 2\sum_{h\in H_{\Gamma}} a(h)\right).
    \end{equation*}
    The degree bound \eqref{eq: F GW L deg} then follows from the identity $$\sum_v (2g(v)-2+|H(v)|) = 2 g - 2 +|L_{\Gamma}|+|D_{\Gamma}|$$ together with the dimension constraint $\sum_{h\in H(v)} a(h) = 3g(v)-3+|H(v)|$, $v\in V_\Gamma$.
\end{proof}

The statement of \cref{thrm:finitley generated n} is not yet in the precise shape of \cref{thm: F GW via Li}\ref{item: thm B i}. By relation \eqref{F g vs correlator} the theorem follows if we can show that \smash{$\langle\!\langle H, \ldots, H \rangle\!\rangle^{Z_n}_{g,m}$} lies in the linear $\mathbb{C}[\ee^n,\dd^n]$-span of
\begin{equation}
    \label{eq: polylog span}
    1, \quad \mr{Li}_{0}\big((-1)^nq\big), \quad \mr{Li}_{-1}\big((-1)^nq\big), \quad \ldots, \quad \mr{Li}_{-N_{g,n} - m}\big((-1)^nq\big)
\end{equation}
with
\begin{equation}
    \label{eq: kmin main}
    N_{g,n} = \left\lfloor \frac{(n-1)(4g-4)}{n}\right\rfloor -1.
\end{equation}
The connection to \cref{thrm:finitley generated n} is established via the relations
\begin{equation}
    \label{eq: polylog vs L}
    \mr{Li}_{0}\big((-1)^nq\big) = \frac{(-1)^n L^n-\ee^n}{\ee^n - (-1)^n \dd^n} , \qquad \mathrm{D} \,\mr{Li}_{-k}\big((-1)^nq\big) = \mr{Li}_{-k-1}\big((-1)^nq\big)
\end{equation}
where the first follows from the defining equation \eqref{eq: def eq L} of $L$.

\begin{proof}[Proof of \cref{thm: F GW via Li}\labelcref{item: thm B i}]
    In the proof of \cref{lem: s t dependence} we already used that all correlators are invariant under permutations of the equivariant variables $\eee_0,\ldots, \eee_{n-1}$ and $\ddd_0,\ldots, \ddd_{n-1}$. Note that a permutation of the first set of variables simultaneously permutes the branches $L_0,\ldots,L_{n-1}$. Thus, as a rational function in the $L_i$, the correlator $$\bigl\langle\!\!\bigl\langle H, \ldots, H \bigr\rangle\!\!\bigr\rangle^{Z_n}_{g,m}$$ is symmetric under permutation of the tuples $(\eee_i,L_i)_{0\leq i < n}$ and permutations of $(\ddd_i)_{0\leq i < n}$. However, since up to a sign its analytic continuation in $q$ is also invariant under the simultaneous exchange of $q\mapsto q^{-1}$ and $\eee_i \leftrightarrow \ddd_i$, $i\in\{0,\ldots,n-1\}$, we conclude that the series must be invariant under permutation of the functions $L_i$ alone. This means, the correlator can be expressed in terms of elementary symmetric polynomials in the $L_i$s. Since $L_i = \zeta^i L$, the only non-vanishing one is $e_n(L_0,\ldots, L_{n-1}) = L^n$ with a similar statement about the $\ee$ and $\dd$-dependence of the Gromov--Witten series (\cref{lem: s t dependence}). We summarise that so far we have proven that
    \begin{equation*}
        \bigl\langle\!\!\bigl\langle H, \ldots, H \bigr\rangle\!\!\bigr\rangle^{Z_n}_{g,m} \in \mathbb{C}(\ee^n,\dd^n) [L^{\pm n}].
    \end{equation*}
    One can furthermore argue that negative powers of $L$ may never occur since due to the $q\rightarrow 0$ expansion
    \begin{equation*}
        L^{-n} = (-1)^n \sum_{d \geq 0} q^d \left(\frac{\dd^{dn}}{\ee^{(d+1)n}} +  \cdots\right)
    \end{equation*}
    they would eventually create denominators in $\ee$ at sufficiently high order in $q$. This is, however, forbidden since, by \cref{lem: s t dependence}, all coefficients of $q^d$ with $d>0$ must be polynomial in $\ee$ and $\dd$. We therefore conclude that the Gromov--Witten series must be polynomial in $L^n$. Since moreover its $L$-degree is bounded, we deduce that the correlator lies in the $\mathbb{C}(\ee^n,\dd^n)$-span of
    \begin{equation*}
        1, \quad L^n, \quad L^{2n}, \quad \ldots, \quad L^{(N_{g,n}+m+1)n}
    \end{equation*}
    with $N_{g,n}$ as in equation \eqref{eq: kmin main}. By relation \eqref{eq: polylog vs L}, this coincides with the subspace spanned by the polylogarithms \eqref{eq: polylog span} concluding our proof of \cref{thm: F GW via Li}\ref{item: thm B i}.
\end{proof}

As a sanity check for our approach we now apply Givental--Teleman reconstruction to $Z_3$ to reprove a formula of Pandharipande--Zinger \cite{PandharipandeNumerative5folds} for the genus-one Gromov--Witten invariants of this fivefold. We already mentioned this formula in the introduction, equation \eqref{eq: PZ formula}, and we remark that a formula for general $n$ can be extracted from our \Cref{thm: closed form}. Such a formula has earlier been proven by Hu \cite{Hu17:g1GWlocPn}.

\begin{example}
    We consider the case $n=3$. Let us apply the Givental--Teleman reconstruction formula to compute
    \begin{equation}
        \mathrm{D} \, \overline{F}{}^{Z_3,T}_{1} = \langle\! \langle H\rangle\! \rangle_{1,1}^{Z_3}\Big|_{\tau_a= \delta_{a,1} \log q} = \sum_{\Gamma \in \mathcal{G}^{[3]}_{1,1}} \mathrm{Cont}(\Gamma).
        \label{sum graphs}
    \end{equation}
    In this case, the stability and the dimension constraints restrict the possibly contributing graphs. Since $D_{i,0} =D_{i,1} =0$, the contribution of graphs with these dilaton leaves vanishes leaving us with only three graphs that may contribute non-trivially. They are shown in \cref{graphs nonzero}.
    
    \begin{figure}[h]
    \centering
        \begin{tikzpicture}[thick, scale= 0.8, transform shape, xshift=-1cm]
        \node[main node] (g1_center) at (0,0) {1};
        \node[right=2.5cm of g1_center] (g1_end) {\Large $H$};
        \draw (g1_center) -- node[above, blue label, pos=0.15] {1} (g1_end);
        \node[caption node] at (1.25, -1.8) {$\Gamma_1$};
        \begin{scope}[xshift=5cm]
            \node[main node] (g2_center) at (0,0) {1};
            \node[right=2.5cm of g2_center] (g2_end) {\Large $H$};
            \node[crossed node] (g2_top) at (0, 2.5) {};
            
            \draw (g2_center) -- node[above, blue label, pos=0.15] {0} (g2_end);
            \draw (g2_center) -- node[left, blue label, pos=0.15] {2} (g2_top);
            \node[caption node] at (1.25, -1.8) {$\Gamma_2$};
        \end{scope}
        \begin{scope}[xshift=11cm]
            \node[main node] (g3_center) at (0,0) {0};
            \node[right=2.5cm of g3_center] (g3_end) {\Large $H$};
            
            \draw (g3_center) -- node[above, blue label, pos=0.15] {0} (g3_end);
            \draw (g3_center) to [out=135, in=225, looseness=8] 
                node[above left, blue label, pos=0.1] {0} 
                node[below left, blue label, pos=0.9] {0} 
                (g3_center);
            \node[caption node] at (1.25, -1.8) {$\Gamma_3$};
            \end{scope}
        \end{tikzpicture}
        \caption{The 3 graphs in $\mathcal{G}_{1,1}^{[2]}$ with a nonzero contribution in the Givental--Teleman reconstruction formula. Cohomological decorations are omitted and the powers of the $\psi$ markings are shown in blue.}
        \label{graphs nonzero}
    \end{figure}
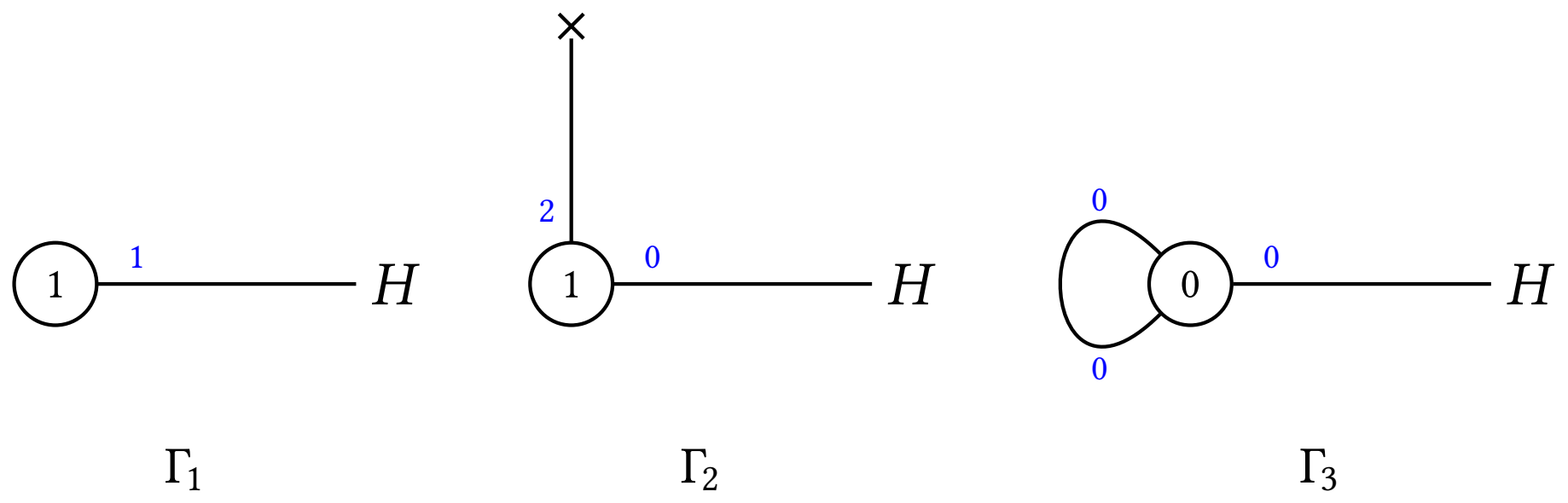
    
    The three graph contributions only require the coefficients $\mathrm{S}_{0i}^{[0]}$ and $\mathrm{S}_{0i}^{[1]}$ of the asymptotic expansion. Solving the Picard–Fuchs equations recursively as in \cref{lemmwwithlogterm} gives
    \begin{align*}
        \mathrm{S}_{0i}^{[0]} &= \frac{c_{0,i}}{L}, \\[4pt]
        \mathrm{S}_{0i}^{[1]} &= \frac{\zeta^{2i}c_{0i}}{36L^5(\ee^3+\dd^3)}\Big(-10L^6-4L^3(\ee^3-\dd^3)+7\ee^3\dd^3+c_{1,i}\Big).
    \end{align*}
    Here, the constants $c_{0i}$ and $c_{1i}$ are chosen to satisfy the initial condition $\mathrm{S}(H^a)|_{q=0} = H^a$. Evaluating the contributions for each graph we find
    \begin{equation}
        \label{eq: D F 1}
        \mathrm{D}\, \overline{F}{}^{Z_3,T}_{1} = -\frac{L^3 - \ee^3 + \dd^3}{8 (\ee^3 + \dd^3)}.
    \end{equation}
    This series includes the constant map contribution
    \[
        \langle H  \rangle_{1,1,0}^{Z_3} = 
        \int_{[\Mbar_{1,1}\times Z_3]_T}H\cdot e(\mathbb{E}^{\vee} \boxtimes T_{Z_3}) = 
        \frac{2\ee^3-\dd^3}{8(\ee^3+\dd^3)}.
        \nonumber
    \]
    If we subtract it from \eqref{eq: D F 1} we obtain
    \[
        \mathrm{D}\, F{}^{Z_3,T}_{1} = - \frac{L^3+\ee^3}{8(\ee^3+\dd^3)} = \frac{1}{8} \, \mr{Li}_0(-q)
    \]
    where for the second equality we use relation \eqref{eq: polylog vs L}. Expanding the polylogarithm yields
    \begin{equation*}
        \int_{[\Mbar_1(Z_3,d)]^{\vir}_T} 1 = \frac{(-1)^d}{8d}
    \end{equation*}
    which matches the formula proven earlier by Pandharipande and Zinger \cite[Lem.~3.1]{PandharipandeNumerative5folds}.
    %
    \label{example:n3}
\end{example}

\section{The anti-diagonal specialisation and proof of \texorpdfstring{\Cref{thm: closed form}}{Theorem C}}
\label{sec: closed form}
This section is dedicated to the proof of the closed formula stated in \Cref{thm: closed form} for the all-genus generating series of equivariant Gromov--Witten invariants of $Z_n$ in fixed degree $d>0$:
\begin{equation}
    \label{eq: fixed-d series}
    \sum_{g \geq 0} u^{2g-2} \int_{[\Mbar_g(Z_n,d)]^{\vir}_{T'}} 1.
\end{equation}
Via virtual localisation \cite{GP97:virtloc}, we will prove a formula for this series in the \emph{anti-diagonal specialisation} of the equivariant parameters:
\begin{equation*}
    \eee_i = - \ddd_i , \qquad i\in\{0,\ldots,n-1\}.
\end{equation*}
We denote by $T'\cong \Gm[n]$ the subtorus of $T$ realising the above weight relations. The key observations which make an evaluation via the localisation formula feasible are:
\begin{itemize}
    \item the fibre of the vector bundle $Z_n\rightarrow \bP^{n-1}$ over any fixed point $P_i$ contains a line with trivial $T'$-action;

    \item in the $T'$-weight decomposition of the tangent space $T_{P_i}Z_n$ at any fixed point, weights always appear in pairs with opposite signs. 
\end{itemize}
\begin{figure}
    \centering
    \begin{tikzpicture}[scale=0.78]
        \node at (-1.5,0.7) {$\boldsymbol{w}_i$:};
        \node at (0,0.7) {\footnotesize$\eee_j - \eee_i$};
        \node at (1.5,0.7) {\footnotesize$\eee_0 - \eee_i$};
        \node at (2.5,0.7) {\footnotesize$\cdots$};
        \node at (3.5,0.7) {\footnotesize$\widehat{\eee_i - \eee_i}$};
        \node at (4.5,0.7) {\footnotesize$\cdots$};
        \node at (5.5,0.7) {\footnotesize$\widehat{\eee_{j} - \eee_i}$};
        \node at (6.5,0.7) {\footnotesize$\cdots$};
        \node at (7.7,0.7) {\footnotesize$\eee_{n-1} - \eee_i$};
        \node at (9.7,0.7) {\footnotesize$-\eee_0 + \eee_i$};
        \node at (10.9,0.7) {\footnotesize$\cdots$};
        \node at (11.5,0.7) {\footnotesize$0$};
        \node at (12.1,0.7) {\footnotesize$\cdots$};
        \node at (13.3,0.7) {\footnotesize$-\eee_{j} + \eee_i$};
        \node at (14.5,0.7) {\footnotesize$\cdots$};
        \node at (15.8,0.7) {\footnotesize$-\eee_{n-1} + \eee_i$};
        \node at (-1.5,-2.7) {$\boldsymbol{w}_j$:};
        \node at (0,-2.7) {\footnotesize$\eee_i - \eee_j$};
        \node at (1.5,-2.7) {\footnotesize$\eee_0 - \eee_j$};
        \node at (2.5,-2.7) {\footnotesize$\cdots$};
        \node at (3.5,-2.7) {\footnotesize$\widehat{\eee_i - \eee_j}$};
        \node at (4.5,-2.7) {\footnotesize$\cdots$};
        \node at (5.5,-2.7) {\footnotesize$\widehat{\eee_{j} - \eee_j}$};
        \node at (6.5,-2.7) {\footnotesize$\cdots$};
        \node at (7.7,-2.7) {\footnotesize$\eee_{n-1} - \eee_j$};
        \node at (9.7,-2.7) {\footnotesize$-\eee_0 + \eee_j$};
        \node at (10.9,-2.7) {\footnotesize$\cdots$};
        \node at (12.1,-2.7) {\footnotesize$-\eee_{i} + \eee_j$};
        \node at (13.3,-2.7) {\footnotesize$\cdots$};
        \node at (13.9,-2.7) {\footnotesize$0$};
        \node at (14.5,-2.7) {\footnotesize$\cdots$};
        \node at (15.8,-2.7) {\footnotesize$-\eee_{n-1} + \eee_j$};
        \node at (0,0) {$\bullet$};
        \draw[thick] (0,0) -- (0,-2);
        \node at (0,-2) {$\bullet$};
        \node[left] at (0,-1) {$\mathbb{P}^{1}$};
        \node[left] at (0,0) {$P_i$};
        \node[left] at (0,-2) {$P_j$};
        \node at (1.5,-1) {$\mathcal{O}_{\mathbb{P}^{1}}(1)$};
        \node at (4.5,-1) {$\cdots$};
        \node at (7.7,-1) {$\mathcal{O}_{\mathbb{P}^{1}}(1)$};
        \node at (9.7,-1) {$\mathcal{O}_{\mathbb{P}^{1}}(-1)$};
        \node at (12.75,-1) {$\cdots$};
        \node at (15.8,-1) {$\mathcal{O}_{\mathbb{P}^{1}}(-1)$};
        \draw [
        thick, 
        decorate, 
        decoration={brace, mirror, raise=0.2cm, amplitude=5pt}
        ] (-0.3,-3.3) -- (8.5,-3.3) 
        node [midway, below, yshift=-0.5cm] {$\mathbb{P}^{n-1}$};
        \draw [
        thick, 
        decorate, 
        decoration={brace, mirror, raise=0.2cm, amplitude=5pt}
        ] (8.9,-3.3) -- (16.9,-3.3) 
        node [midway, below, yshift=-0.5cm] {fibres of $Z_n \rightarrow \mathbb{P}^{n-1}$};
    \end{tikzpicture}
    \caption{The normal bundle of a torus-preserved $\mathbb{P}^{1}$ inside $Z_n$ connecting two fixed points $P_i$ and $P_j$ together with the tangent $T'$-weights $\boldsymbol{w}_i=(w^k_i)_{0\leq k \leq 2n-2}$ at the fixed points.}
    \label{fig: line wts}
\end{figure}
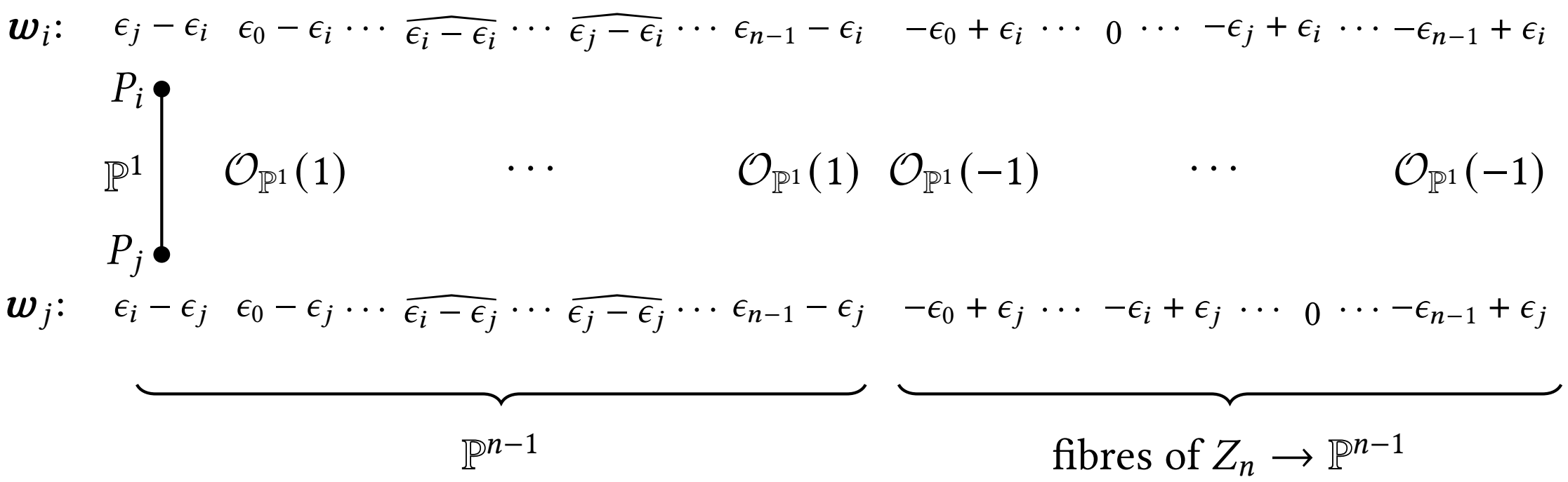
This is illustrated in \Cref{fig: line wts} which displays the $T'$-weights over two fixed points in $Z_n$ joined by a $T'$-preserved line. The first property ensures that if we evaluate the fixed-degree series \eqref{eq: fixed-d series} via virtual localisation, the only non-vanishing contributions may come from $T'$-fixed stable maps whose domain features a rational component giving a $d$-fold ramified covering of a torus-preserved line in $Z_n$ fully ramified over the torus fixed points it connects. See \cite[Sec.~2]{FP00:HodgeIntGW} or \cite[Sec.~3]{PandharipandeNumerative5folds} for a similar calculation. This means all remaining domain components must be contracted to either of the two fixed points. With knowledge of the tangent $T'$ weights as displayed in \Cref{fig: line wts}, one can show that the fixed-degree series \eqref{eq: fixed-d series} gets the following contribution from a torus-preserved line joining two fixed points $P_i$ and $P_j$, $i\neq j$, of $Z_n$:
\begin{itemize}
    \item The degree-$d$ ramified cover of the torus-preserved line contributes a factors
    \begin{equation*}
        (-1)^{(d+1)n} \frac{1}{d} \frac{1}{\prod_{k\neq i}(\eee_k-\eee_i)} \, \frac{1}{\prod_{k\neq j}(\eee_k-\eee_j)}.
    \end{equation*}

    \item Contracted domain components to $P_i$ contribute
    \begin{equation}
        \label{eq: fxd pt contrib}
        \sum_{g \geq 0} u^{2g-1} \int_{\Mbar_{g,1}} \frac{\prod_{k=1}^{2n-1} \Lambda^{\vee}_g(w_i^{k})}{\frac{\eee_j - \eee_i}{d} - \psi_1}, \qquad \quad \Lambda^{\vee}_g(\epsilon) = \sum_{l\geq 0} \lambda_l (-1)^{l} \epsilon^{g-l}
    \end{equation}
    with a suitable unstable genus-zero term.

    \item The contribution from contracted components to $P_j$ is analogous.
\end{itemize}
The last two contributions can be simplified using Mumford's relation \cite[Sec.~5]{Mum83:Towards}
\begin{equation*}
    \Lambda^{\vee}_g(\epsilon) \cdot \Lambda^{\vee}_g(-\epsilon) = (-1)^g \epsilon^{2g}
\end{equation*}
and the special evaluation
\begin{equation*}
    \sum_{g \geq 0} u^{2g-1} \int_{\Mbar_{g,1}} \psi_1^{2g-2} \lambda_g = \frac{1}{2 \sin \frac{u}{2}}
\end{equation*}
proven by Faber and Pandharipande \cite[Thm.~2]{FP00:HodgeIntGW}. Inserting the tangent weights displayed in \Cref{fig: line wts} into expression \eqref{eq: fxd pt contrib} thus leads to the following simplified contribution:
\begin{equation*}
    \sum_{g \geq 0} u^{2g-1} \! \int_{\Mbar_{g,1}} \!\!\frac{\prod_{k=1}^{2n-1} \Lambda^{\vee}_g(w_i^{k})}{\frac{\eee_j - \eee_i}{d} - \psi_1} = \begin{cases}
        {\displaystyle\frac{1}{2 \sin \big(\frac{du}{2} \prod_{l \neq i,j} (\eee_l-\eee_i)\big)} \, \prod_{k\neq i}(\eee_k-\eee_i)} & \text{if }n\text{ is even,} \\[1.7em]
        {\displaystyle\frac{1}{2 \sinh \big(\frac{du}{2}  \prod_{l \neq i,j} (\eee_l-\eee_i)\big)} \, \prod_{k\neq i}(\eee_k-\eee_i)} & \text{if }n\text{ is odd.}
    \end{cases}
\end{equation*}
Combining the contributions of the ramified cover and the contracted components and summing over all torus-preserved lines between two fixed points in $Z_n$ thus yields
\begin{equation*}
    \sum_{g \geq 0} u^{2g-2} \int_{[\Mbar_g(Z_n,d)]^{\vir}_{T'}} 1 = \frac{1}{d} \sum_{0\leq i<j <n} ~ \frac{1}{\prod_{k\in\{i,j\}} 2 \sin \big(\frac{du}{2}  \prod_{l\neq i,j} (\eee_l - \eee_k)\big)}
\end{equation*}
when $n$ is even and
\begin{equation}
    \label{eq: closed form n odd}
    \sum_{g \geq 0} u^{2g-2} \int_{[\Mbar_g(Z_n,d)]^{\vir}_{T'}} 1 = \frac{(-1)^{d+1}}{d} \sum_{0\leq i<j <n} ~ \frac{1}{\prod_{k\in\{i,j\}} 2 \sinh \big(\frac{du}{2} \prod_{l\neq i,j} (\eee_l - \eee_k)\big)}
\end{equation}
for odd $n$. This closes our proof of \Cref{thm: closed form}. We remark that the simplifications arising from the specific choice of torus weights are a direct generalisation to arbitrary $n$ of the simplifications exploited earlier by Faber--Pandharipande in their proof of the $n=2$ case \cite[Sec.~2]{FP00:HodgeIntGW} and Pandharipande--Zinger's calculation for $n=3$ \cite[Sec.~3]{PandharipandeNumerative5folds}.

\begingroup
\renewcommand*{\bibfont}{\normalfont\small}
\printbibliography
\endgroup
\end{document}